\documentstyle{article}

\newcommand{\bigzerou}{%
\smash{\lower1.7ex\hbox{\bg 0}}}
\newtheorem{theorem}{Theorem}
\newtheorem{prop}{Proposition}
\newtheorem{defi}{Definition}
\newtheorem{cor}{Corollary}
\newtheorem{conj}{Conjecture}
\newtheorem{Rem}{Remark}

\newcommand{\ba}{\begin{eqnarray}}
\newcommand{\ea}{\end{eqnarray}}
\newcommand{\no}{\nonumber}

\def\d{{\partial}}
\newcommand{\mapright}[1]{%
\smash{\mathop{%
\hbox to 1.0cm{\rightarrowfill}}\limits^{#1}}}
\newcommand{\mapleft}[1]{%
\smash{\mathop{%
\hbox to 1.3cm{\leftarrowfill}}\limits^{#1}}}

\begin{document}
\title{
\begin{flushright}
  \begin{minipage}[b]{5em}
    \normalsize
    ${}$      \\
  \end{minipage}
\end{flushright}
{\bf Coordinate Change of Gauss-Manin System and Generalized Mirror 
Transformation\\}}
\author{Masao Jinzenji\\
\\
\it Division of Mathematics, Graduate School of Science \\
\it Hokkaido University \\
\it  Kita-ku, Sapporo, 060-0810, Japan\\
{\it e-mail address: jin@math.sci.hokudai.ac.jp}}
\maketitle
\begin{abstract}
In this paper, we explicitly derive the generalized mirror transformation of 
quantum cohomology of general type projective hypersurfaces, proposed 
in our previous article, as an effect of coordinate change of 
the virtual Gauss-Manin system.
\end{abstract}
\section{Introduction}
This paper is the completion of our previous work \cite{gene}, 
\cite{vir} on the quantum cohomology 
of general type hypersurfaces.
To be more precise, we construct an 
algorithm to compute K\"ahler Gromov-Witten invariants (structure constants of
K\"ahler sub-ring of small quantum cohomology ring) of degree $k$ 
hypersurface in $CP^{N-1}$
with $k>N$ (we denote it by $M_{N}^{k}$) for rational curves of 
arbitrary degree 
only by using Givental's ODE \cite{giv}
\begin{equation}
\biggl((\d_{x})^{N-1}-k\cdot e^{x}\cdot (k\d_{x}+k-1)(k\d_{x}+k-2)
\cdots(k\d_{x}+1)\biggr)w(x)=0,
\label{ifun}
\end{equation}
as the starting point.\footnote{ From now on, we omit the word ``K\"ahler'' 
and limit ourselves to consider Gromov-Witten invariants of $M_{N}^{k}$ with 
insertion of ${\cal O}_{e^{j}}$, where $e$ is the K\"ahler class of 
$M_{N}^{k}$.}
Of course, there are already several literatures
\cite{gath}, \cite{blly} and \cite{givc} on the quantum cohomology ring 
of general type projective hypersurfaces. But in \cite{gath} 
and \cite{blly}, connection with mirror 
computation or (\ref{ifun}) seems to be implicit. As for \cite{givc}, 
direct connection with hypergeometric class 
or with (\ref{ifun}) is suggested, but their construction is a little abstract 
and not practical for explicit 
prediction. From this point of view, our construction presented 
in this paper is explicit and can be used to  
compute any genus $0$ Gromov-Witten invariants of $M_{N}^{k}$ 
directly from (\ref{ifun}).

 Our construction is 
fundamentally based on the line of thoughts exposed 
in the course of our previous papers \cite{cj},
\cite{fano},\cite{gene},\cite{vir} and \cite{gm}. 
Main ingredients of our construction consist of 
virtual structure constants and generalized mirror 
transformation, which were introduced and partially 
constructed in \cite{gene}. 
The virtual structure constants have their origin in our work \cite{cj} with
A.Collino. In \cite{cj}, 
A.Collino and myself constructed recursive formulas that express the structure 
constants of small quantum cohomology ring of 
$M_{N}^{k}\;\;(N-k)\geq 2$ in terms of the ones of $M_{N+1}^{k}$ up to 
degree $5$ rational curves (explicit form of 
recursive formula for rational curves of arbitrary degree
was conjectured in \cite{fano} and was proved in \cite{gm}). 
The virtual structure constants are defined 
as rational number obtained from iterative use of 
these recursive formulas to the region $(N-k\leq 1)$.
In \cite{gene}, we conjectured that the virtual structure 
constants are deeply related to (\ref{ifun}) 
and that they correspond to three point functions of 
the B-model in the mirror conjecture. These conjecture 
were proved to be true in \cite{gm}. In other words, 
we have constructed an algorithm to compute the virtual
structure constants for arbitrary $N$ and $k$ only 
by using (\ref{ifun}). These results correspond to the first 
part of the mirror computation of the small quantum 
cohomology of general type projective hypersurface 
$M_{N}^{k}\;\;(k>N)$. 

The second part of the mirror computation is 
the generalized mirror transformation, 
that correspond to the formula to convert 
the virtual structure constants into the actual structure constants 
of the small quantum cohomology ring of $M_{N}^{k}$ 
in our context. Note that in the Calabi-Yau case, this 
process of mirror computation is realized by 
the coordinate change of the B-model deformation parameter into
flat coordinate of the A-model, or K\"ahler 
deformation parameter. The notion of the generalized mirror 
transformation for $M_{N}^{k}\;\;(k>N)$ was first 
introduced in \cite{gene} and explicit form of the 
generalized mirror transformation was determined up 
to degree $3$ rational curves by comparing the value 
of the virtual structure constants with the numerical 
results of some three point Gromov-Witten invariants 
of $M_{N}^{k}$. In \cite{vir}, we determined the form of 
the generalized mirror transformation up to degree 
$5$ rational curves by using curious coincidence 
between the terms appearing in the results of \cite{gene}
and many-point Gromov-Witten invariants obtained from 
application of the associativity equation \cite{km} and 
the K\"ahler equation to some simple combination of 
the virtual structure constants. This idea turned out to
be very effective, but unfortunately, we found that 
such coincidence was not complete and that some non-trivial
modification was needed as degree of rational 
curves grows. Therefore, our search for the generalized 
mirror transformation temporarily stopped.

Next break-through comes from application of 
the Gauss-Manin system to the quantum cohomology ring of $M_{N}^{k}$
\cite{gm}. In \cite{gm}, we proposed an idea of the 
virtual Gauss-Manin system, that has the virtual structure 
constants as the matrix elements of the Gauss-Manin connection. 
The virtual Gauss-Manin system, when applied to 
the Calabi-Yau hypersurfaces $M_{k}^{k}$, 
is effective not only in the B-model computation, but also in finding 
out the Jacobian between the B-model deformation parameter 
and the flat coordinate of the A-model, as was 
suggested in the last section of \cite{gm}. 
Therefore, we searched for a way to extend the construction of the 
mirror transformation of $M_{k}^{k}$ via coordinate 
change of the virtual Gauss-Manin system to the case of 
$M_{N}^{k}\;\;(k>N)$. In the sequel, a hint of the answer 
was given by the theory of Iritani \cite{iri}, which 
is a concrete exposition of the idea of Givental and Coates \cite{givc}. 
According to Iritani, we have to introduce not only 
the K\"ahler deformation parameter $x^{1}$, but also the parameter
$x^{j}\;(j=2,\cdots,N-2)$ that corresponds to $e^{j}$, 
i.e., $j$ times wedge product of the K\"ahler class $e$.
With this idea, we extend the virtual 
Gauss-Manin system by constructing the virtual Gauss-Manin 
connection corresponding to the deformation by $x^{j}$. 
With these preparation, we can read off the matrix 
elements of the Jacobian $\frac{\d t^{j}}{\d x^{i}}$ between 
$x^{j}$ and the A-model deformation parameter
$t^{j}$ from the matrix elements of the virtual 
Gauss-Manin connection. Moreover, we can integrate out the 
Jacobian and obtain the rule of coordinate change 
$x^{i}=x^{i}(t^{1},\cdots,t^{N-2})$. The results so obtained 
force us to introduce many point Gromov-Witten 
invariants obtained from application of the associativity equation
and the modified K\"ahler equation, that was derived from 
Iritani's framework, to the matrix elements 
of Gauss-Manin connection associated with deformation of 
$t^{1}$. After some straightforward computations, we 
can naturally reproduce the results up to degree $5$ 
rational curves given in \cite{vir}. Of course, our 
construction given in this paper can be applied to 
predict arbitrary genus $0$ Gromov-Witten invariants of 
$M_{N}^{k}\;\;(k>N)$. Indeed, we computed some three point Gromov-Witten 
invariants of $M_{13}^{14}$ for degree $6$ rational 
curves by using our construction, and checked coincidence 
with the numerical results obtained from fixed point 
computation by Kontsevich \cite{kont}. We also expect that 
our construction will be explained from compactification of 
moduli space of rational curves in $CP^{N-1}$ \cite{gene},\cite{bert}.

This paper is organized as follows. In Section 2, 
we introduce the notation, the definition and the results 
obtained in our previous works. Especially, the definition 
of the virtual structure constants, the virtual 
Gauss-Manin system and the virtual Gromov-Witten 
invariants are given, and their relation to Givental`s ODE 
(\ref{ifun}) is discussed. We also mention some 
conjectures on the form of the generalized mirror transformation 
 given in \cite{gene}. In Section 3, we expose our 
algorithm to construct the generalized mirror transformation for rational 
curves of arbitrary degree from the virtual Gauss-Manin system. 
Next, we derive our previous results in \cite{vir}.
Note that our presentation in this section is 
slightly different from the outline given here, because 
we want to respect priority of the work of Iritani. 
Of course, the result of computation does not change. 
In Section 4, we briefly review the results of 
Iritani and explain how our algorithm in Section 3 naturally 
follows from his framework. 
  \\
\\  
{\bf Acknowledgement} We would like to thank 
Prof. T.Eguchi and  Prof. M.Noumi for discussions at an early stage of 
this work. We also thank Prof. Y.Tonegawa and  Hokkaido University Computer 
Center for assistance on numerical computation using fixed point 
theorem. Finally, we especially thank Dr. H.Iritani for suggestions that 
led me to a key idea to finish this work. 
\section{Overview of Our Previous Results}
\subsection{Definitions}
In this subsection, we introduce the quantum K\"ahler sub-ring 
of the quantum cohomology ring of a degree $k$ hypersurface in
$CP^{N-1}$.
Let $M_{N}^{k}$ be a hypersurface of degree $k$ in $CP^{N-1}$.
 We denote by $QH^{*}_{e}(M_{N}^{k})$ the 
subring of the quantum cohomology ring $QH^{*}(M_{N}^{k})$
generated by ${\cal O}_{e}$ induced from the K\"ahler form $e$ 
(or, equivalently the intersection $H\cap M_{N}^{k}$ between a hyperplane
class $H$ of $CP^{N-1}$ and $M_{N}^{k}$).
Additive basis of $QH_{e}^{*}(M_{N}^{k})$ is given by 
${\cal O}_{e^{j}}\;\;(j=0,1,\cdots,N-2)$, which is induced from 
$e^{j}\in H^{j,j}(M_{N}^{k})$.
 The multiplication rule of $QH^{*}_{e}(M_{N}^{k})$ 
is determined by the Gromov-Witten invariant of genus $0$ 
$\langle {\cal O}_{e}{\cal O}_{{e}^{N-2-m}}
{\cal O}_{{e}^{m-1-(k-N)d}}\rangle_{d,M_{N}^{k}}$ and
it is given as follows:
\begin{eqnarray}
 L_{m}^{N,k,d} &:=&\frac{1}{k}\langle {\cal O}_{e}{\cal O}_{{e}^{N-2-m}}
{\cal O}_{{e}^{m-1-(k-N)d}}\rangle_{d},\no\\
\no\\
{\cal O}_{e}\cdot 1&=&{\cal O}_{e},\nonumber\\
{\cal O}_{e}\cdot{\cal O}_{{e}^{N-2-m}}&=&{\cal O}_{{e}^{N-1-m}}+
\sum_{d=1}^{\infty}L_{m}^{N,k,d}q^{d}{\cal O}_{{e}^{N-1-m+(k-N)d}},\no\\
q&:=&\exp(t),
\label{gm}
\end{eqnarray}
where the subscript $d$ counts the degree of the rational curves
measured by $e$. Therefore,  $q=\exp(t)$ is the degree counting 
parameter. 
\begin{defi}
We call $L_{n}^{N,k,d}$ the structure constant of weighted degree $d$.
\end{defi}
Since $M_{N}^{k}$ is a complex $(N-2)$ dimensional manifold, we see that
a structure constant $L_{m}^{N,k,d}$
is non-zero only if the following condition is satisfied:
\begin{eqnarray}
&& 1\leq N-2-m\leq N-2, 1\leq m-1+(N-k)d\leq N-2,\no\\
&\Longleftrightarrow &max\{0,2-(N-k)d\}\leq m \leq min\{N-3,N-1-(N-k)d\}.
\label{sel}
\end{eqnarray}
We rewrite (\ref{sel}) into 
\begin{eqnarray}
(N-k\geq 2) &\Longrightarrow& 0\leq m \leq (N-1)-(N-k)d\no\\
(N-k=1,d=1)&\Longrightarrow& 1\leq m \leq N-3\no\\
(N-k=1,d\geq2)&\Longrightarrow& 0\leq m \leq N-1-(N-k)d\no\\
(N-k\leq 0)&\Longrightarrow& 2+(k-N)d\leq m \leq N-3.
\label{flasel}
\end{eqnarray}
From (\ref{flasel}), we easily see that the number of the non-zero
structure constants $L_{m}^{N,k,d}$ is finite except for the case of $N=k$.
Moreover, if $N\geq 2k$, the non-zero structure constants come only from
the $d=1$ part and the non-vanishing $L_{m}^{N,k,1}$  
is determined by $k$ and  
independent of $N$. 
\subsection{Virtual Structure Constants and Givental's ODE}
In this subsection, we introduce the virtual structure constants 
$\tilde{L}^{N,k,d}_{m}$ which is non-zero only if $0\leq m\leq N-1+(k-N)d$.
The original definition of $\tilde{L}^{N,k,d}_{m}$ in \cite{gm} 
is given by the initial condition:
\begin{equation}
\sum_{m=0}^{k-1}\tilde{L}_{m}^{N,k,1}w^{m}=
k\cdot\prod_{j=1}^{k-1}(jw+(k-j)),\;\;(N-k\geq 2),
\label{vi}
\end{equation}
and the recursive formulas that describe $\tilde{L}^{N,k,d}_{m}$ as 
a weighted homogeneous polynomial in $\tilde{L}^{N+1,k,d'}_{n}\;\;(d'\leq d)$
of degree $d$. See \cite{fano} 
for the explicit form of the recursive formulas. In \cite{gm}, we showed that the virtual structure constants are directly connected with the Givental's
ODE:     
\begin{equation}
\biggl((\d_{x})^{N-1}-k\cdot e^{x}\cdot (k\d_{x}+k-1)(k\d_{x}+k-2)\cdots(k\d_{x}+1)\biggr)w(x)=0,
\label{fun}
\end{equation}
for arbitrary $N$ and $k$ via the virtual Gauss-Manin system defined
as follows: 
\begin{defi}
We call the following rank 1 ODE for vector valued function 
$\tilde{\psi}_{m}(x)$
\begin{eqnarray}
\partial_{x}\tilde{\psi}_{N-2-m}(x)&=&\tilde{\psi}_{N-1-m}(x)+
\sum_{d=1}^{\infty}
\exp(dx)\cdot \tilde{L}_{m}^{N,k,d}\cdot\tilde{\psi}_{N-1-m-(N-k)d}(x).
\label{gm1}
\end{eqnarray}
the virtual Gauss-Manin system associated with the quantum K\"ahler 
sub-ring of $M_{N}^{k}$, where $m$ runs through $0\leq m \leq N-2$ 
if $N-k \geq 1$, $0\leq m \leq N-1$ if $N-k=0$, and $m\in {\bf Z}$ if 
$N-k<0$.  
\end{defi}
Here, we restate the main result in \cite{gm}.
\begin{theorem}
We can derive the following identity from the virtual Gauss-Manin system 
(\ref{gm1}).
\begin{equation}
\tilde{\psi}_{N-1}(x)=\biggl((\partial_{x})^{N-1}-k\cdot e^{x}\cdot
(k\partial_{x}+k-1)\cdots
(k\partial_{x}+2)\cdot(k\partial_{x}+1)\biggr)\partial_{x}^{\beta}
\tilde{\psi}_{-\beta}(x)
\label{hyper} 
\end{equation}
where $\beta=0$ if $N-k\geq 1$, $\beta=1$ if $N-k=0$, and $\beta=\infty$ 
if $N-k<0$.
\end{theorem}
We can also compute $\tilde{L}^{N,k,d}_{m}$ only by using the above theorem,
 and this process is the analogue of the B-model computation of the 
Calabi-Yau case.
\begin{cor}
The virtual structure constants $\tilde{L}_{n}^{N,k,d}$ are fully 
reconstructed from the identity (\ref{hyper}).
As a result, we can compute all the virtual structure constants by using the 
relation:
\begin{eqnarray}
&&\sum_{n=0}^{k-1}\tilde{L}_{n}^{N,k,1}w^{n}=
k\cdot\prod_{j=1}^{k-1}(jw+(k-j)), \no\\
&&\sum_{m=0}^{N-1+(k-N)d}\tilde{L}_{m}^{N,k,d}z^{m}=\no\\
&&\sum_{l=2}^{d}(-1)^{l}\sum_{0=i_{0}<\cdots<i_{l}=d}
\sum_{j_{l}=0}^{N-1+(k-N)d}\cdots\sum_{j_{2}=0}^{j_{3}}\sum_{j_{1}=0}^{j_{2}}
\prod_{n=1}^{l}\biggl((\frac{i_{n-1}+(d-i_{n-1})z}{d})^{j_{n}-j_{n-1}}
\cdot \tilde{L}_{j_{n}+(N-k)i_{n-1}}^{N,k,i_{n}-i_{n-1}}\biggr).\no\\
\label{cinic}
\end{eqnarray}
\end{cor}
Note that the condition for non-vanishing $\tilde{L}_{m}^{N,k,d}$:
$0\leq m\leq N-1+(k-N)d$, implies that we have infinite number 
of $\tilde{L}_{m}^{N,k,d}$'s in the $k>N$ case. 
In \cite{gene}, we proposed a conjecture of the formula 
to convert the virtual structure constants $\tilde{L}_{m}^{N,k,d}$
into the real structure constants $L_{m}^{N,k,d}$. To restate the conjecture 
in \cite{gene}, we introduce some combinatorial definitions. 
\begin{defi}
Let $P_{d}$ be the set of partitions of positive integer $d$:
\begin{equation}
P_{d}=\{\sigma_{d}=(d_{1},d_{2},\cdots,d_{l(\sigma_{d})})\;\;|\;\;
1\leq d_{1}\leq d_{2}\leq\cdots\leq d_{l(\sigma_{d})}\;\;,\;\;\
\sum_{j=1}^{l(\sigma_{d})}d_{j}=d\;\;,\;\;d_{j}\in{\bf Z}\}.
\label{part} 
\end{equation}
From now on, we denote a partition $\sigma_{d}$ by 
$d_{1}+d_{2}+\cdots+d_{l(\sigma_{d})}$. In (\ref{part}), we denote 
the length of the partition $\sigma_{d}$ by $l(\sigma_{d})$.
If $d=0$, we define $P_{0}$ as the set that consists of one trivial 
partition $0$ with length $0$. 
\end{defi}
\begin{defi}
Let $S(\sigma_d)$ be a rational number 
associated with the partition $\sigma_{d}\in
P_{d}$, which is defined by the following generating function:
\begin{equation}
\sum_{d=0}^{\infty}
(\sum_{\sigma_d\in P_{d}}S(\sigma_d)
\prod_{j=1}^{l(\sigma_{d})}a_{d_{j}})z^{d}
:=\exp(\sum_{j=1}^{\infty}a_{j}z^{j}).
\label{not1}
\end{equation}
\end{defi}
With these preparation, we restate the main conjecture in \cite{gene}
on the structure of the generalized mirror transformation:
\begin{conj}
The generalized mirror transformation takes the form
\begin{eqnarray}
L^{N,k,d}_{n}&=&\sum_{m=0}^{d-1}\sum_{\sigma_{m}\in P_{m}}
(-1)^{l(\sigma_{m})}\cdot d^{l(\sigma_{m})}\cdot S(\sigma_{m})
\cdot\prod_{i=1}^{l(\sigma_{m})}
\biggl(\frac{\tilde{L}_{1+(k-N)d_{i}}^{N,k,d_{i}}}{d_{i}}
\biggr)\cdot G_{d-m}^{N,k,d}(n;\sigma_{m}),\no\\
\label{gene}
\end{eqnarray}
where $G_{d-m}^{N,k,d}(n;\sigma_{m})$ is a polynomial of 
$\tilde{L}_{n}^{N,k,d}$  with weighted degree $d$. 
\end{conj}
Of course, we have to determine the polynomial $G_{d-m}^{N,k,d}(n;\sigma_{m})$
to predict the real structure constants $L_{n}^{N,k,d}$. In \cite{vir}, 
we have determined $G_{d-m}^{N,k,d}(n;\sigma_{m})$ for arbitrary $N$ and $k$ 
up to $d=3$ case, and $G_{d-m}^{k-1,k,d}(n;\sigma_{m})$ up to $d=5$ 
case. In these cases, we only need $\tilde{L}_{n}^{N,k,d}$'s  that satisfy 
$1+(k-N)d\leq n \leq N-2$. Therefore, we only have to use finite 
number of $\tilde{L}_{n}^{N,k,d}$'s in the $k >N$ case in spite of the fact 
that we have infinite number of the virtual structure constants. This 
observation plays an important role in our construction given in the 
next section. Moreover, the above conjecture naturally follows from 
the construction.
 We will show the explicit algorithm to determine the unknown 
polynomial $G_{d-m}^{N,k,d}(n;\sigma_{m})$ for arbitrary degree $d$ 
in the next section.
\subsection{Virtual Gromov-Witten Invariants}
We introduce here the virtual Gromov-Witten invariants to make
correspondence of the results of this article with the ones 
in \cite{vir}. Precisely speaking, these quantities do not play
an essential role in the main result of this article, 
but they are convenient for describing our results 
of computation. 
\begin{defi}
The virtual Gromov-Witten invariant
$v(\prod_{j=1}^{n}{\cal O}_{e^{a_{j}}})_{d}$ on $M_{N}^{k}$ is the 
rational number that satisfy the condition:\\
(i) initial condition
\begin{eqnarray}
&&v({\cal O}_{e^{a}}{\cal O}_{e^{b}}{\cal O}_{e^{c}})_{0}=
k\cdot\delta_{a+b+c,N-2},\no\\
&&v(\prod_{j=1}^{n}{\cal O}_{e^{a_{j}}})_{0}=0,
\;\;(n\neq 3),\no\\
&&\frac{1}{k}v({\cal O}_{e^{N-2-n}}{\cal O}_{e^{n-1-(k-N)d}}{\cal O}_{e})_{d}=
\tilde{L}_{n}^{N,k,d}-\tilde{L}_{1+(k-N)d}^{N,k,d},\;\;(d\geq 1),
\end{eqnarray}
(ii) flat metric condition
\begin{eqnarray}
&&v({\cal O}_{e^{0}}{\cal O}_{e^{a}}{\cal O}_{e^{b}})_{0}
=k\cdot\delta_{a+b,N-2},\no\\
&&v({\cal O}_{e^{0}}\prod_{j=1}^{n}{\cal O}_{e^{a_{j}}})_{d}=0,
\;\;(d\geq 1,\;\;\mbox{or} \;\;\;d=0,\;\; n\neq 2),
\end{eqnarray}
(iii) topological selection rule\\
\begin{equation}
v(\prod_{j=1}^{n}{\cal O}_{e^{a_{j}}})_{d}\neq0
\Longrightarrow (N-5)+(N-k)d=\sum_{j=1}^{n}(a_{j}-1),
\end{equation}
(iv) K\"ahler equation\\
\begin{equation}
v({\cal O}_{e}\prod_{j=1}^{n}{\cal O}_{e^{a_{j}}})_{d}=
d\cdot v(\prod_{j=1}^{n}{\cal O}_{e^{a_{j}}})_{d},
\label{simk}
\end{equation}
(v) associativity equation\\
\begin{eqnarray}
&&\sum_{d_{1}=0}^{d}\sum_{\{\alpha_{*}\}\coprod\{\beta_{*}\}=\{n_{*}\}}\sum_{i=0}^{N-2}v({\cal O}_{e^{a}}{\cal O}_{e^{b}}
(\prod_{\alpha_{j}\in\{\alpha_{*}\}}{\cal O}_{e^{\alpha_{j}}}){\cal O}_{e^{i}})_{d_{1}}
v({\cal O}_{e^{N-2-i}}  
(\prod_{\beta_{j}\in\{\beta_{*}\}}{\cal O}_{e^{\beta_{j}}}){\cal O}_{e^{c}}
{\cal O}_{e^{d}})_{d-d_{1}}\no\\
&&=\sum_{d_{1}=0}^{d}\sum_{\{\alpha_{*}\}\coprod\{\beta_{*}\}=\{n_{*}\}}\sum_{i=0}^{N-2}
v({\cal O}_{e^{a}}{\cal O}_{e^{c}}
(\prod_{\alpha_{j}\in\{\alpha_{*}\}}{\cal O}_{e^{\alpha_{j}}}){\cal O}_{e^{i}})_{d_{1}}
v({\cal O}_{e^{N-2-i}}  
(\prod_{\beta_{j}\in\{\beta_{*}\}}{\cal O}_{e^{\beta_{j}}}){\cal O}_{e^{b}}
{\cal O}_{e^{d}})_{d-d_{1}},\no\\
&&(a+b+c+d+\sum_{j=1}^{m}(n_{j}-1)=N-2+(N-k)d).
\end{eqnarray}
\end{defi}
Next, we introduce the notation that was heavily used in \cite{vir},
\begin{defi}
\begin{equation}
V_{d-m}^{N,k,d}(n;d_{1}+d_{2}+\cdots+d_{l(\sigma_{m})})
:=\frac{1}{k\cdot (d-m)^{l(\sigma_{m})-1}}
v({\cal O}_{e^{N-2-n}}{\cal O}_{e^{n-1-(k-N)d}}\prod_{j=1}^{l(\sigma_{m})}
{\cal O}_{e^{1+(k-N)d_{j}}})_{d-m}.
\end{equation}
\end{defi}
In \cite{vir}, we proposed some conjectures on the relation 
between $V_{d-m}^{N,k,d}(n;\sigma_{m})$ and 
$G_{d-m}^{N,k,d}(n;\sigma_{m})$ in (\ref{gene}).
\begin{conj}
If $l(\sigma_{m})\leq 1$ or $d-m=1$, $G_{d-m}^{N,k,d}(n;\sigma_{m})$
is given by $V_{d-m}^{N,k,d}(n;\sigma_{m})$.
\end{conj}
In the same way as Conjecture 1, we can see that this conjecture is 
natural and consistent with the construction given in the next section.
As the last part of this subsection, we write down a recursive formula 
for $V_{d-m}^{k-1,k,d}(n;\sigma_{m})$, that follows directly
from the definition 
of the virtual structure constants, for later use:  
 \begin{prop}
\begin{eqnarray}
&&V_{d-m}^{k-1,k,d}(n;d_{1}+d_{2}+\cdots+d_{l(\sigma_{m})})\no\\
&&=V_{d-m}^{k-1,k,d-1}(n;d_{1}+d_{2}+\cdots+d_{l(\sigma_{m})-1}+
(d_{l(\sigma_{m})}-1))\no\\
&&+V_{d-m}^{k-1,k,d-d_{l(\sigma_{m})}}(n-d_{l(\sigma_{m})}
;d_{1}+d_{2}+\cdots+d_{l(\sigma_{m})-1})\no\\
&&-V_{d-m}^{k-1,k,d-d_{l(\sigma_{m})}}(d+1
;d_{1}+d_{2}+\cdots+d_{l(\sigma_{m})-1})\no\\
&&+\sum_{j=1}^{d-m-1}\sum_{A\coprod B}(\frac{d-m-j}{d-m})^{p}
(\frac{j}{d-m})^{l(\sigma_{m})-p-1}V_{d-m-j}^{k-1,k,d-m-j+d_{A}+
d_{l(\sigma_{m})}-1}(n;d_{a_{1}}+\cdots+d_{a_{p}}+(d_{l(\sigma_{m})}-1))
\no\\
&&\times V_{j}^{k-1,k,j+d_{B}}
(n+j-d+m-d_{A}-d_{l(\sigma_{m})};d_{b_{1}}+\cdots+d_{b_{l(\sigma_{m})-1-p}})
\no\\
&&-\sum_{j=1}^{d-m-1}\sum_{A\coprod B}(\frac{d-m-j}{d-m})^{p}
(\frac{j}{d-m})^{l(\sigma_{m})-p-1}V_{d-m-j}^{k-1,k,d}
(n;d_{a_{1}}+\cdots+d_{a_{p}}+(m+j-d_{A}))
\no\\
&&\times V_{j}^{k-1,k,j+d_{B}}
(m+1+j-d_{A};d_{b_{1}}+\cdots+d_{b_{l(\sigma_{m})-1-p}}),
\label{vgw}
\end{eqnarray}
where $\sum_{A\coprod B}$ means the summation on all the way of 
separating the set $\{1,2,\cdots,l(\sigma_{m})-1\}$ into two disjoint sets 
$A=\{a_{1},\cdots,a_{p}\}$ and $B=\{b_{1},\cdots,b_{l(\sigma_{m})-1-p}\}$,
and $d_{A}$ (resp. $d_{B}$) is an integer 
$\sum_{j=1}^{p}d_{a_{j}}$ (resp. $\sum_{j=1}^{l(\sigma_{m})-1-p}d_{b_{j}}$).
\end{prop}
\begin{Rem}
In applying (\ref{vgw}), we don't need to arrange integers $d_{j}$ in 
ascending order, because the virtual Gromov-Witten invariants  are 
invariant under permutation of insertion of operators. 
\end{Rem}
\section{Algorithmic Derivation of the Generalized Mirror Transformation}
In this section, we write down the algorithm to compute the generalized
mirror transformation and reproduce the results obtained in \cite{vir} 
for $M_{k-1}^{k}$ up to $d=5$ case. We explain theoretical background 
of our algorithm, following the result of Iritani \cite{iri}, in the next 
section. 

As the first step, we truncate the virtual Gauss-Manin system (\ref{gm1}) 
into the form:
\begin{eqnarray}
\partial_{x}\psi_{N-2-m}(x)&=&\psi_{N-1-m}(x)+
\sum_{d=1}^{[\frac{m-1}{k-N}]}
\exp(dx)\cdot \tilde{L}_{m}^{N,k,d}\cdot\psi_{N-1-m+(k-N)d}(x).\no\\
&&(0\leq m\leq N-2)
\label{tgm}
\end{eqnarray}
This truncation means that we throw away all the $\tilde{L}^{N,k,d}_{m}$'s
except for the ones that satisfy $(1+(k-N)d\leq m\leq N-2)$.
Moreover, $\tilde{\psi}_{j}(x)\;(j\in{\bf Z})$ is replaced by 
$\psi_{j}(x)\;(j=0,\cdots N-2)$. In other words, we throw away  
$\tilde{\psi}_{j}(x)$'s, that are exotic as the usual Gauss-Manin  
system associated with the quantum cohomology ring of $M_{N}^{k}$.
For example, 
we write down the form of truncated virtual Gauss-Manin system of $M_{8}^{9}$:
 \begin{eqnarray}
\d_{x}\psi_{0}(x)&=&\psi_{1}(x)+\alpha e^{x}\psi_{2}(x)+
\eta e^{2x}\psi_{3}(x)+\varphi e^{3x}\psi_{4}(x)
+\pi e^{4x}\psi_{5}(x)+\epsilon e^{5x}\psi_{6}(x)\no\\
\d_{x}\psi_{1}(x)&=&\psi_{2}(x)+\beta e^{x}\psi_{3}(x)+
\xi e^{2x}\psi_{4}(x)+\kappa e^{3x}\psi_{5}(x)+\pi e^{4x}\psi_{6}(x)\no\\
\d_{x}\psi_{2}(x)&=&\psi_{3}(x)+\gamma e^{x}\psi_{4}(x)+
\xi e^{2x}\psi_{5}(x)+\varphi e^{3x}\psi_{6}(x)\no\\
\d_{x}\psi_{3}(x)&=&\psi_{4}(x)+\beta e^{x}\psi_{5}(x)+
\eta e^{2x}\psi_{6}(x)\no\\
\d_{x}\psi_{4}(x)&=&\psi_{5}(x)+\alpha e^{x}\psi_{6}(x)\no\\
\d_{x}\psi_{5}(x)&=&\psi_{6}(x)\no\\
\d_{x}\psi_{6}(x)&=&0.
\label{89}
\end{eqnarray}
Since the next step is rather complicated, we explain it by using the 
above example. First, we eliminate $\psi_{6}(x)$ from the 5th line of 
(\ref{89}) by using the 6th line,
\begin{eqnarray}
\d_{x}\psi_{4}(x)&=&\psi_{5}(x)+\alpha e^{x}\d_{x}\psi_{5}(x)
\no\\
&=& (1+\alpha e^{x}\d_{x})\psi_{5}(x),
\label{dmo1}
\end{eqnarray}
and obtain, 
\begin{eqnarray}
\psi_{5}(x)&=&(1+\alpha e^{x}\d_{x})^{-1}\d_{x}\psi_{4}(x),\no\\
\psi_{6}(x)&=&\d_{x}(1+\alpha e^{x}\d_{x})^{-1}\d_{x}\psi_{4}(x).
\label{dmo2}
\end{eqnarray}
In this way, we can eliminate $\psi_{5}(x)$ and $\psi_{6}(x)$ in (\ref{89}).
We can also eliminate $\psi_{4}(x)$ by plugging (\ref{dmo2}) into 
4th line of (\ref{89}),
\begin{eqnarray}
\d_{x}\psi_{3}(x)=\biggl(1+\beta e^{x}(1+\alpha e^{x}\d_{x})^{-1}\d_{x}
+\eta e^{2x}\d_{x}(1+\alpha e^{x}\d_{x})^{-1}\d_{x}\biggr)\psi_{4}(x),
\label{dmo3}
\end{eqnarray} 
and obtain,
\begin{eqnarray}
&&\psi_{4}(x)=\biggl(1+\beta e^{x}(1+\alpha e^{x}\d_{x})^{-1}\d_{x}
+\eta e^{2x}\d_{x}(1+\alpha e^{x}\d_{x})^{-1}\d_{x}\biggr)^{-1}
\d_{x}\psi_{3}(x),\no\\
&&\psi_{5}(x)=(1+\alpha e^{x}\d_{x})^{-1}\d_{x}
\biggl(1+\beta e^{x}(1+\alpha e^{x}\d_{x})^{-1}\d_{x}
+\eta e^{2x}\d_{x}(1+\alpha e^{x}\d_{x})^{-1}\d_{x}\biggr)^{-1}
\d_{x}\psi_{3}(x),\no\\
&&\psi_{6}(x)=\d_{x}(1+\alpha e^{x}\d_{x})^{-1}\d_{x}
\biggl(1+\beta e^{x}(1+\alpha e^{x}\d_{x})^{-1}\d_{x}
+\eta e^{2x}\d_{x}(1+\alpha e^{x}\d_{x})^{-1}\d_{x}\biggr)^{-1}
\d_{x}\psi_{3}(x).
\label{dmo4}
\end{eqnarray}
By continuing the elimination process, we can obtain 
the following relation
between $\psi_{0}(x)$ and $\psi_{1}(x)$:
\begin{equation}
\psi_{1}(x)=F(e^{x},\d_{x})\psi_{0}(x).
\label{flat}
\end{equation}
Then we dictate the following relation that characterize the flat 
coordinate $t$,
\begin{equation}
\d_{t}=F(e^{x},\d_{x}),\;\;\mbox{i.e.,}\;\;\;\; \d_{t}\psi_{0}(x)=\psi_{1}(x).
\label{flat2}
\end{equation} 
Explicitly, $\d_{t}$ in our example is given as follows: 
\begin{eqnarray}
\d_{t}&=&\d_{x}-\alpha e^{x}\d_{x}^{2}+e^{2x}((\alpha\beta-\eta)\d_{x}^{3}
+\alpha^{2}(\d_{x}+1)\d_{x}^{2})
-e^{3x}(\varphi \d_{x}^{4}-
\eta\beta(\d_{x}+1)\d_{x}^{3}-\eta\gamma\d_{x}^{4}+(\alpha\beta\gamma-
\alpha\xi)\d_{x}^{4}\no\\
&&+\alpha\beta^{2}(\d_{x}+1)\d_{x}^{3}
+\alpha(\alpha\beta-\eta)((\d_{x}+2)\d_{x}^{3}+(\d_{x}+1)^{2}\d_{x}^{2})+
\alpha^{3}(\d_{x}+2)(\d_{x}+1)\d_{x}^{2})+\cdots.
\label{gmt}
\end{eqnarray} 
Using (\ref{gmt}), we can rewrite (\ref{89}) into the form,
\begin{eqnarray}
\d_{t}\psi_{0}(x)&=&\psi_{1}(x)\no\\
\d_{t}\psi_{1}(x)&=&\psi_{2}(x)+(\beta-\alpha)(e^{x}-
\alpha e^{2x}+(\alpha\beta-\eta+
2\alpha^{2})e^{3x})\psi_{3}(x)\no\\
&&+((\xi-\eta-\alpha(\gamma-\alpha))(e^{2x}-2\alpha e^{3x})+
((\alpha\beta-\eta+\alpha^{2})(\beta-\alpha)+ 
(\alpha\beta-\eta)(\gamma-\alpha))e^{3x})\psi_{4}(x)\no\\
&&+(\kappa-\varphi-\alpha(\xi-\eta)-\eta(\beta-\alpha)+
\alpha^{2}(\gamma-\alpha))e^{3x}\psi_{5}(x)\no\\
\d_{t}\psi_{2}(x)&=&\psi_{3}(x)+(\gamma-\alpha)(e^{x}-
\alpha e^{2x}+(\alpha\beta-\eta+2\alpha^{2})e^{3x})\psi_{4}(x)\no\\
&&+((\xi-\eta-\alpha(\gamma-\alpha))(e^{2x}-2\alpha e^{3x})+
(2(\alpha\beta-\eta)+\alpha^{2})(\gamma-\alpha))e^{3x})\psi_{5}(x)\no\\
\d_{t}\psi_{3}(x)&=&\psi_{4}(x)+(\beta-\alpha)(e^{x}-
\alpha e^{2x}+(\alpha\beta-\eta+
2\alpha^{2})e^{3x})\psi_{5}(x)+(\beta-\alpha)(\alpha\beta-\eta+\alpha^{2})
\psi_{6}(x)\no\\
\d_{t}\psi_{4}(x)&=&\psi_{5}(x)\no\\
\d_{t}\psi_{5}(x)&=&\psi_{6}(x)\no\\
\d_{t}\psi_{6}(x)&=&0
\label{1st}
\end{eqnarray}
Though these computations are rather hard because of non-commutativity 
between $e^{x}$ and $\d_{x}$, we can manage them with the aid of Maple
package: Ore-algebra.  
In general, we obtain the following equations as the result of this step,
\begin{eqnarray}
\d_{t}\psi_{0}(x)&=&\psi_{1}(x),\no\\
\d_{t}\psi_{N-2-m}(x)&=&\psi_{N-1-m}(x)+
\sum_{d=1}^{[\frac{m-2}{k-N}]}f^{N,k,d}_{m}(e^{x})\cdot\psi_{N-1-m+d(k-N)}(x),
\;\;(1\leq m \leq N-3),\no\\
\d_{t}\psi_{N-2}(x)&=&0,
\end{eqnarray}
where $f^{N,k,d}_{m}(e^{x})$ is the power series in $e^{x}$ whose 
lowest power is more than $d$. 
\begin{defi}
We denote by $\frac{1}{k}w({\cal O}_{e^{N-2-m}}{\cal O}_{e^{m-1-(k-N)d}}
{\cal O}_{e})_{d}$
the coefficient of $e^{dx}$ in $f^{N,k,d}_{m}(e^{x})$.
\end{defi}
In our example, they are given as follows:
\begin{eqnarray}
&&\frac{1}{k}w({\cal O}_{e}{\cal O}_{e}{\cal O}_{e^{3}})_{1}=\beta-\alpha,
\;\;
\frac{1}{k}w({\cal O}_{e}{\cal O}_{e^{2}}{\cal O}_{e^{2}})_{1}=\gamma-\alpha,
\no\\
&&\frac{1}{k}w({\cal O}_{e}{\cal O}_{e}{\cal O}_{e^{2}})_{2}=
\xi-\eta-\alpha(\gamma-\alpha),\no\\
&&\frac{1}{k}w({\cal O}_{e}{\cal O}_{e}{\cal O}_{e})_{3}=
\kappa-\varphi-\alpha(\xi-\eta)-\eta(\beta-\alpha)+\alpha^{2}(\gamma-\alpha).
\label{wex1}
\end{eqnarray}
Explicitly, we can write down 
$\frac{1}{k}w({\cal O}_{e^{N-2-m}}{\cal O}_{e^{m-1-(k-N)d}}{\cal O}_{e})_{d}$
in terms of virtual Gromov-Witten invariants up to $d=3$ cases for 
arbitrary $N$ and $k$.
\begin{eqnarray}
\frac{1}{k}w({\cal O}_{e^{N-2-m}}{\cal O}_{e^{m-1-(k-N)}}{\cal O}_{e})_{1}
&=&V_{1}^{N,k,1}(n;0),\no\\
\frac{1}{k}w({\cal O}_{e^{N-2-m}}{\cal O}_{e^{m-1-2(k-N)}}{\cal O}_{e})_{2}
&=&V_{2}^{N,k,2}(n;0)-\tilde{L}^{N,k,1}_{1+(k-N)}V_{1}^{N,k,2}(n;1),\no\\
\frac{1}{k}w({\cal O}_{e^{N-2-m}}{\cal O}_{e^{m-1-3(k-N)}}{\cal O}_{e})_{3}
&=&V_{3}^{N,k,3}(n;0)-\tilde{L}^{N,k,1}_{1+(k-N)}V_{2}^{N,k,3}(n;1)
-\tilde{L}^{N,k,2}_{1+2(k-N)}V_{1}^{N,k,3}(n;2)\no\\
&&+(\tilde{L}^{N,k,1}_{1+(k-N)})^{2}V_{1}^{N,k,3}(n;1+1).
\label{wex2}
\end{eqnarray}
We also show here the explicit results for $d=4,5$ rational curves 
in the case of $N=k-1$:
\begin{eqnarray}
\frac{1}{k}w({\cal O}_{e^{k-3-m}}{\cal O}_{e^{m-5}}{\cal O}_{e})_{4}
&=&V_{4}^{k-1,k,4}(n;0)-\tilde{L}^{k-1,k,1}_{2}V_{3}^{k-1,k,4}(n;1)
-\tilde{L}^{k-1,k,2}_{3}V_{2}^{k-1,k,4}(n;2)\no\\
&&-\tilde{L}^{k-1,k,3}_{4}V_{1}^{k-1,k,4}(n;3)\no\\
&&+(\tilde{L}^{k-1,k,1}_{2})^{2}\biggl(
V_{2}^{k-1,k,3}(n;1)+V_{2}^{k-1,k,3}(n-1;1)
-V_{2}^{k-1,k,2}(5;0)+A(n)\biggr)\no\\
&&+2\tilde{L}^{k-1,k,1}_{2}\tilde{L}^{k-1,k,2}_{3}
V_{1}^{k-1,k,4}(n;1+2)\no\\
&&-(\tilde{L}^{k-1,k,1}_{2})^{3}V_{1}^{k-1,k,4}(n;1+1+1),\no\\
\frac{1}{k}w({\cal O}_{e^{k-3-m}}{\cal O}_{e^{m-6}}{\cal O}_{e})_{5}
&=&V_{5}^{k-1,k,5}(n;0)-\tilde{L}^{k-1,k,1}_{2}V_{4}^{k-1,k,5}(n;1)
-\tilde{L}^{k-1,k,2}_{3}V_{3}^{k-1,k,5}(n;2)\no\\
&&-\tilde{L}^{k-1,k,3}_{4}V_{4}^{k-1,k,5}(n;3)
-\tilde{L}^{k-1,k,4}_{5}V_{1}^{k-1,k,5}(n;4)\no\\
&&+(\tilde{L}_{2}^{k-1,k,1})^{2}\biggl(
V_{3}^{k-1,k,4}(n;1)+V_{3}^{k-1,k,4}(n-1;1) 
-V_{3}^{k-1,k,3}(6;0)\no\\
&&+B(n)+C(n)-V_{1}^{k-1,k,1}(3;0)\cdot
(2hi_{1}(n)+2hi_{2}(n)+hi_{3}(n))\biggr)\no\\
&&+2\tilde{L}_{2}^{k-1,k,1}\tilde{L}_{3}^{k-1,k,2}
\biggl(V_{2}^{k-1,k,4}(n;2)+V_{2}^{k-1,k,4}(n-1;2)
-V_{2}^{k-1,k,2}(6;0)\no\\
&&+2hi_{1}(n)+\frac{3}{2}hi_{2}(n)+hi_{3}(n)-\frac{1}{2}hi_{4}(n)
\biggr)\no\\
&&-(\tilde{L}_{2}^{k-1,k,1})^{3}
\biggl(V_{2}^{k-1,k,3}(n;1)+2V_{2}^{k-1,k,3}(n-1;1)+V_{2}^{k-1,k,3}(n-2;1)\no\\
&&-2V_{2}^{k-1,k,2}(5;0)-V_{2}^{k-1,k,3}(6;1)
+D(n)+3hi_{1}(n)+3hi_{2}(n)+hi_{3}(n)\biggr)\no\\
&&-3(\tilde{L}_{2}^{k-1,k,1})^{2}\tilde{L}_{3}^{k-1,k,2}
V_{1}^{k-1,k,5}(n;1+1+2)\no\\
&&+(\tilde{L}_{2}^{k-1,k,1})^{4}V_{1}^{k-1,k,5}(n;1+1+1+1),
\label{wex}
\end{eqnarray}
where
\begin{eqnarray}
A(n):&=&V_{1}^{k-1,k,2}(n;1)\cdot V_{1}^{k-1,k,1}(n-3;0)
+V_{1}^{k-1,k,1}(n;0)\cdot
V_{1}^{k-1,k,2}(n-2;1)\no\\
&&-V_{1}^{k-1,k,4}(n;1+2)\cdot
V_{1}^{k-1,k,1}(3;0)
-V_{1}^{k-1,k,4}(n;3)\cdot
V_{1}^{k-1,k,1}(4;0),\no\\
B(n):&=&V_{2}^{k-1,k,3}(n;1)\cdot V_{1}^{k-1,1,1}(n-4;0)
+V_{1}^{k-1,k,1}(n;0)\cdot V_{2}^{k-1,k,3}(n-2;1)\no\\
&&-\bigl(V_{2}^{k-1,k,4}(n;2)+V_{2}^{k-1,k,4}(n-1;2)
-V_{2}^{k-1,k,2}(6;0)\bigr)\cdot V_{1}^{k-1,k,1}(3;0)\no\\
&&-V_{1}^{k-1,k,5}(n;4)\cdot V_{2}^{k-1,k,2}(5;0),\no\\
C(n)&:=&V_{1}^{k-1,k,2}(n;1)\cdot V_{2}^{k-1,k,2}(n-3;0)+
V_{2}^{k-1,k,2}(n;0)\cdot V_{1}^{k-1,k,2}(n-3;1)\no\\
&&-V_{1}^{k-1,k,5}(n;1+3)\cdot V_{2}^{k-1,k,2}(4;0)
-V_{2}^{k-1,k,5}(n;3)
\cdot V_{1}^{k-1,k,1}(4;0),\no\\
D(n)&:=&V_{1}^{k-1,k,2}(n;1)\cdot V_{1}^{k-1,k,1}(n-3;0)
+V_{1}^{k-1,k,1}(n;0)\cdot
V_{1}^{k-1,k,2}(n-2;1)\no\\
&&-V_{1}^{k-1,k,4}(n;1+2)\cdot
V_{1}^{k-1,k,1}(3;0)
-V_{1}^{k-1,k,4}(n;3)\cdot
V_{1}^{k-1,k,1}(4;0)\no\\
&&+V_{1}^{k-1,k,2}(n-1;1)\cdot V_{1}^{k-1,k,1}(n-4;0)
+V_{1}^{k-1,k,1}(n-1;0)\cdot
V_{1}^{k-1,k,2}(n-3;1)\no\\
&&-V_{1}^{k-1,k,4}(n-1;1+2)\cdot V_{1}^{k-1,k,1}(3;0)
-V_{1}^{k-1,k,4}(n-1;3)\cdot V_{1}^{k-1,k,1}(4;0).
\label{hidden}
\end{eqnarray}
In (\ref{wex}), $hi_{j}(n)$ is a degree $2$ homogeneous polynomial of 
$\tilde{L}^{k-1,k,1}_{m}$ satisfying $hi_{j}(6)=hi_{j}(7)=0$ and is
given by,  
\begin{eqnarray}
hi_{1}(n)&=&\tilde{L}_{n}^{k-1,k,1}\tilde{L}_{n-4}^{k-1,k,1}-
\tilde{L}_{3}^{k-1,k,1}(\tilde{L}_{n}^{k-1,k,1}+\tilde{L}_{n-1}^{k-1,k,1}+
\tilde{L}_{n-2}^{k-1,k,1}+\tilde{L}_{n-3}^{k-1,k,1}+\tilde{L}_{n-4}^{k-1,k,1})
\no\\
&&+\tilde{L}_{2}^{k-1,k,1}(\tilde{L}_{n-1}^{k-1,k,1}+\tilde{L}_{n-2}^{k-1,k,1}+\tilde{L}_{n-3}^{k-1,k,1})\no\\
&&-\bigl(\tilde{L}_{6}^{k-1,k,1}\tilde{L}_{2}^{k-1,k,1}-
\tilde{L}_{3}^{k-1,k,1}(\tilde{L}_{6}^{k-1,k,1}+\tilde{L}_{5}^{k-1,k,1}+
\tilde{L}_{4}^{k-1,k,1}+\tilde{L}_{3}^{k-1,k,1}+\tilde{L}_{2}^{k-1,k,1})\no\\
&&+\tilde{L}_{2}^{k-1,k,1}(\tilde{L}_{5}^{k-1,k,1}+
\tilde{L}_{4}^{k-1,k,1}+\tilde{L}_{3}^{k-1,k,1})\bigr),\no\\
hi_{2}(n)&=&\tilde{L}_{n}^{k-1,k,1}\tilde{L}_{n-3}^{k-1,k,1}
+\tilde{L}_{n-1}^{k-1,k,1}\tilde{L}_{n-4}^{k-1,k,1}\no\\
&&-\tilde{L}_{4}^{k-1,k,1}(\tilde{L}_{n}^{k-1,k,1}+\tilde{L}_{n-1}^{k-1,k,1}+
\tilde{L}_{n-2}^{k-1,k,1}+\tilde{L}_{n-3}^{k-1,k,1}+\tilde{L}_{n-4}^{k-1,k,1})
+\tilde{L}_{2}^{k-1,k,1}\tilde{L}_{n-2}^{k-1,k,1}\no\\
&&-\bigl(\tilde{L}_{6}^{k-1,k,1}\tilde{L}_{3}^{k-1,k,1}
+\tilde{L}_{5}^{k-1,k,1}\tilde{L}_{2}^{k-1,k,1}\no\\
&&-\tilde{L}_{4}^{k-1,k,1}(\tilde{L}_{6}^{k-1,k,1}+\tilde{L}_{5}^{k-1,k,1}+
\tilde{L}_{4}^{k-1,k,1}+\tilde{L}_{3}^{k-1,k,1}+\tilde{L}_{2}^{k-1,k,1})
+\tilde{L}_{2}^{k-1,k,1}\tilde{L}_{4}^{k-1,k,1}\bigr),\no\\
hi_{3}(n)&=&\tilde{L}_{n}^{k-1,k,1}\tilde{L}_{n-2}^{k-1,k,1}
+\tilde{L}_{n-1}^{k-1,k,1}\tilde{L}_{n-3}^{k-1,k,1}+
\tilde{L}_{n-2}^{k-1,k,1}\tilde{L}_{n-4}^{k-1,k,1}\no\\
&&-\tilde{L}_{5}^{k-1,k,1}(\tilde{L}_{n}^{k-1,k,1}+\tilde{L}_{n-1}^{k-1,k,1}+
\tilde{L}_{n-2}^{k-1,k,1}+\tilde{L}_{n-3}^{k-1,k,1}+\tilde{L}_{n-4}^{k-1,k,1})
-\tilde{L}_{2}^{k-1,k,1}\tilde{L}_{n-2}^{k-1,k,1}\no\\
&&-\bigl(\tilde{L}_{6}^{k-1,k,1}\tilde{L}_{4}^{k-1,k,1}
+\tilde{L}_{5}^{k-1,k,1}\tilde{L}_{3}^{k-1,k,1}+
\tilde{L}_{4}^{k-1,k,1}\tilde{L}_{2}^{k-1,k,1}\no\\
&&-\tilde{L}_{5}^{k-1,k,1}(\tilde{L}_{6}^{k-1,k,1}+\tilde{L}_{5}^{k-1,k,1}+
\tilde{L}_{4}^{k-1,k,1}+\tilde{L}_{3}^{k-1,k,1}+\tilde{L}_{2}^{k-1,k,1})
-\tilde{L}_{2}^{k-1,k,1}\tilde{L}_{4}^{k-1,k,1}\bigr),\no\\
hi_{4}(n)&=&\tilde{L}_{n-1}^{k-1,k,1}\tilde{L}_{n-3}^{k-1,k,1}
-\tilde{L}_{4}^{k-1,k,1}(\tilde{L}_{n-1}^{k-1,k,1}+
\tilde{L}_{n-2}^{k-1,k,1}+\tilde{L}_{n-3}^{k-1,k,1})
+\tilde{L}_{3}^{k-1,k,1}\tilde{L}_{n-2}^{k-1,k,1}\no\\
&&-\bigl(\tilde{L}_{5}^{k-1,k,1}\tilde{L}_{3}^{k-1,k,1}
-\tilde{L}_{4}^{k-1,k,1}(\tilde{L}_{5}^{k-1,k,1}+
\tilde{L}_{4}^{k-1,k,1}+\tilde{L}_{3}^{k-1,k,1})
+\tilde{L}_{3}^{k-1,k,1}\tilde{L}_{4}^{k-1,k,1}\bigr).\no\\
\end{eqnarray}
For later use, we also write down the explicit formula of
some  $V_{d-m}^{k-1,k,d}(n;\sigma_{m})$'s \footnote{In \cite{vir},
we made a combinatorial mistake in evaluating 
$V_{2}^{k-1,k,5}(n;1+1+1)$, and the terms including $hi_{j}(n)$ are 
different from the fourth equality in (\ref{virc}). Therefore,
(\ref{virc}) is the correct answer.}:  
\begin{eqnarray}
V_{2}^{k-1,k,4}(n;1+1)&=&V_{2}^{k-1,k,3}(n;1)+V_{2}^{k-1,k,3}(n-1;1)
-V_{2}^{k-1,k,2}(5;0)+\frac{1}{2}A(n)\no\\
V_{3}^{k-1,k,5}(n;1+1)&=&V_{3}^{k-1,k,4}(n;1)+V_{3}^{k-1,k,4}(n-1;1) 
-V_{3}^{k-1,k,3}(6;0),\no\\
&&+\frac{2}{3}B(n)+\frac{1}{3}C(n)-\frac{1}{3}V_{1}^{k-1,k,1}(3;0)\cdot
(2hi_{1}(n)+hi_{2}(n)+hi_{3}(n)-hi_{4}(n)),\no\\
V_{2}^{k-1,k,5}(n;1+2)&=&V_{2}^{k-1,k,4}(n;2)+V_{2}^{k-1,k,4}(n-1;2)
-V_{2}^{k-1,k,2}(6;0)\no\\
&&+hi_{1}(n)+\frac{1}{2}hi_{2}(n)+\frac{1}{2}hi_{3}(n)-\frac{1}{2}hi_{4}(n),
\no\\
V_{2}^{k-1,k,5}(n;1+1+1)&=&V_{2}^{k-1,k,3}(n;1)+2V_{2}^{k-1,k,3}(n-1;1)+V_{2}^{k-1,k,3}(n-2;1)\no\\
&&-2V_{2}^{k-1,k,2}(5;0)-V_{2}^{k-1,k,3}(6;1)\no\\
&&+\frac{1}{2}D(n)+\frac{1}{4}(4hi_{1}(n)+4hi_{2}(n)+hi_{3}(n)+hi_{4}(n)).
\label{virc} 
\end{eqnarray}
We can obtain (\ref{virc}) by direct application of (\ref{vgw}). 
It is quite interesting that the structure of the terms proportional 
to $\prod_{j=1}^{l(\sigma_{m})}\tilde{L}^{N,k,d}_{1+(k-n)d_{j}}$ 
in $\frac{1}{k}
w({\cal O}_{e^{N-2-n}}{\cal O}_{e^{n-1-(k-N)d}}{\cal O}_{e})_{d}$
resembles the one of $V_{d-m}^{N,k,d}(n;\sigma_{m})$.
We would like the reader to compare 
the corresponding $G_{d-m}^{k-1,k,d}(n;\sigma_{m})$'s in the conjecture 
(\ref{gene}), that were determined in \cite{vir}:
\begin{eqnarray}
G_{2}^{k-1,k,4}(n;1+1)&=&V_{2}^{k-1,k,3}(n;1)+V_{2}^{k-1,k,3}(n-1;1)
-V_{2}^{k-1,k,2}(5;0)+\frac{3}{4}A(n)\no\\
G_{3}^{k-1,k,5}(n;1+1)&=&V_{3}^{k-1,k,4}(n;1)+V_{3}^{k-1,k,4}(n-1;1) 
-V_{3}^{k-1,k,3}(6;0),\no\\
&&+\frac{4}{5}B(n)+\frac{3}{5}C(n)-\frac{1}{5}V_{1}^{k-1,k,1}(3;0)\cdot
(6hi_{1}(n)+5hi_{2}(n)+3hi_{3}(n)-hi_{4}(n)),\no\\
G_{2}^{k-1,k,5}(n;1+2)&=&V_{2}^{k-1,k,4}(n;2)+V_{2}^{k-1,k,4}(n-1;2)
-V_{2}^{k-1,k,2}(6;0)\no\\
&&+\frac{1}{5}(8hi_{1}(n)+5hi_{2}(n)+4hi_{3}(n)-3hi_{4}(n)),
\no\\
G_{2}^{k-1,k,5}(n;1+1+1)&=&V_{2}^{k-1,k,3}(n;1)+2V_{2}^{k-1,k,3}(n-1;1)+V_{2}^{k-1,k,3}(n-2;1)\no\\
&&-2V_{2}^{k-1,k,2}(5;0)-V_{2}^{k-1,k,3}(6;1)\no\\
&&+\frac{4}{5}D(n)+
\frac{1}{25}(46hi_{1}(n)+46hi_{2}(n)+16hi_{3}(n)-2hi_{4}(n)).
\label{virg} 
\end{eqnarray}
As was suggested in \cite{vir}, we can see that the coefficients of 
$A(n)$,$B(n)$, $C(n)$ and $D(n)$ in (\ref{virc}) change followingly 
in (\ref{virg}):
\begin{equation}
\frac{1}{2}\rightarrow \frac{3}{4},\;\;\; 
\frac{2}{3}\rightarrow \frac{4}{5},\;\;\;
\frac{1}{3}\rightarrow \frac{3}{5},\;\;\; 
\frac{1}{2}\rightarrow \frac{4}{5}.
\label{change}  
\end{equation} 
The terms that consist of $hi_{j}(n)$ change in a more complicated manner.
These changes seem to imply that $G_{d-m}^{k-1,k,d}(n;\sigma_{m})$ 
does not satisfy the simple K\"ahler equation given in (\ref{simk}). 
At this stage, we introduce a key definition to resolve these puzzles 
observed in \cite{vir}.
\begin{defi}
Let $w({\cal O}_{e^{n_1}}{\cal O}_{e^{n_2}}\cdots {\cal O}_{e^{n_m}})$
be the rational number obtained from applying the associativity equation
and the modified K\"ahler equation:
\begin{eqnarray}
&&w({\cal O}_{e}{\cal O}_{e^{n_{1}}}{\cal O}_{e^{n_{2}}}\cdots
{\cal O}_{e^{n_{m}}})_{d}\nonumber\\
&&=d\cdot 
w({\cal O}_{e^{n_{1}}}{\cal O}_{e^{n_{2}}}\cdots{\cal O}_{e^{n_{m}}})_{d}-
\sum_{f=1}^{d-1}\tilde{L}_{1+(k-N)f}^{N,k,f}\cdot
w({\cal O}_{e^{1+(k-N)f}}{\cal O}_{e^{n_{1}}}{\cal O}_{e^{n_{2}}}\cdots
{\cal O}_{e^{n_{m}}})_{d-f},
\label{mok}
\end{eqnarray}
to the initial condition $w({\cal O}_{e^{N-2-m}}{\cal O}_{e^{m-1-(k-N)d}}
{\cal O}_{e})_{d}$.
\end{defi}
\begin{prop}
\begin{equation}
\frac{1}{k}w({\cal O}_{e^{N-2-n}}{\cal O}_{e^{n-1-(k-N)d}}{\cal O}_{e}
\prod_{i=1}^{l(\sigma_{d-1})}{\cal O}_{e^{1+(k-N)d_{i}}})_{1}
=V_{1}^{N,k,d}(n;d_{1}+d_{2}+\cdots+d_{l(\sigma_{d-1})}).
\label{d1}
\end{equation}
\end{prop}
{\it proof)}
In $d=1$ case, the modified K\"ahler equation (\ref{mok}) reduces to
(\ref{simk}). Moreover, the initial condition and the associativity 
equation are also the same as the ones for the virtual Gromov-Witten 
invariants. Therefore, the assertion of the proposition follows. $\Box$

With these set up, our main result of this article is given as follows:
\begin{conj}
\begin{eqnarray}
&&\langle{\cal O}_{e^{N-2-m}}{\cal O}_{e^{m-1-(k-N)d}}{\cal O}_{e}
\rangle_{d}\no\\
&&=\sum_{g=0}^{d-1}(-1)^{l(\sigma_{g})}\sum_{\sigma_{g}\in P_{g}}
S(\sigma_{g})
\biggl(\prod_{i=1}^{l(\sigma_{g})}
\frac{\tilde{L}^{N,k,d_{i}}_{1+(k-N)d_{i}}}{d_{i}}\biggr)
w({\cal O}_{e^{N-2-m}}{\cal O}_{e^{m-1-(k-N)d}}{\cal O}_{e}
\prod_{i=1}^{l(\sigma_{g})}{\cal O}_{e^{1+(k-N)d_{i}}})_{d-g}.
\no\\
\label{main}
\end{eqnarray}
\end{conj}
\begin{Rem}
Under the assumption of (\ref{main}), (\ref{mok}) leads us to the standard 
K\"ahler equation for real Gromov-Witten invariants:
\begin{eqnarray}
&&\langle{\cal O}_{e}
{\cal O}_{e^{N-2-m}}{\cal O}_{e^{m-1-(k-N)d}}{\cal O}_{e}
\rangle_{d}\no\\
&&=\sum_{g=0}^{d-1}(-1)^{l(\sigma_{g})}\sum_{\sigma_{g}\in P_{g}}
S(\sigma_{g})
\biggl(\prod_{i=1}^{l(\sigma_{g})}
\frac{\tilde{L}^{N,k,d_{i}}_{1+(k-N)d_{i}}}{d_{i}}\biggr)
w({\cal O}_{e}{\cal O}_{e^{N-2-m}}{\cal O}_{e^{m-1-(k-N)d}}{\cal O}_{e}
\prod_{i=1}^{l(\sigma_{g})}{\cal O}_{e^{1+(k-N)d_{i}}})_{d-g}
\no\\
&&=d\cdot\langle{\cal O}_{e^{N-2-m}}{\cal O}_{e^{m-1-(k-N)d}}{\cal O}_{e}
\rangle_{d}.
\label{rk}
\end{eqnarray}
In order to derive (\ref{rk}), we introduce another representation of 
$\sigma_{d-g}\in P_{d-g}$:
\begin{equation}
\sigma_{g}=d_{1}+\cdots+d_{l(\sigma_{g})}=\sum_{j=1}^{\infty}
m_{j}\cdot j.
\label{mulr}
\end{equation}
Then the combinatorial factor $S(\sigma_{g})
\prod_{i=1}^{l(\sigma_{g})}\frac{1}{d_{i}}$  is rewritten as follows:
\begin{equation}
S(\sigma_{g})
\prod_{i=1}^{l(\sigma_{g})}\frac{1}{d_{i}}=
\prod_{j=1}^{\infty}\frac{1}{(m_{j})!\cdot j^{m_{j}}}.
\label{com}
\end{equation}
In these set up, what we have to show is,
\begin{eqnarray}
(d-g)\prod_{j=1}^{\infty}\frac{1}{(m_{j})!\cdot j^{m_{j}}}
+\sum_{i\;(m_{i}>0)} m_{i}\cdot i
\prod_{j=1}^{\infty}\frac{1}{(m_{j})!\cdot j^{m_{j}}}=
d\prod_{j=1}^{\infty}\frac{1}{(m_{j})!\cdot j^{m_{j}}},
\label{comde}
\end{eqnarray}
but this identity directly follows from 
$\sum_{i\;(m_{i}>0)}m_{i}\cdot i=g$.
\end{Rem}
In the remaining part of this section,
we explicitly reproduce the results in \cite{vir} by using Conjecture 3. 
First, we derive the results of arbitrary $N$ and $k$ $(N-k<0)$
up to $d=3$ case. In the $d=1$ case, (\ref{main}) merely says,
\begin{equation}
\frac{1}{k}\langle{\cal O}_{e^{N-2-n}}{\cal O}_{e^{n-1-(k-N)}}{\cal O}_{e}
\rangle_{1}
=\frac{1}{k}w({\cal O}_{e^{N-2-n}}{\cal O}_{e^{n-1-(k-N)}}{\cal O}_{e})_{1}
=V_{1}^{N,k,1}(n;0)=\tilde{L}_{n}^{N.k.1}-\tilde{L}_{1+(k-N)}^{N,k,1}.
\label{r1}
\end{equation}
In the $d=2$ case, we can derive 
\begin{eqnarray}
&&\frac{1}{k}\langle{\cal O}_{e^{N-2-n}}{\cal O}_{e^{n-1-2(k-N)}}{\cal O}_{e}
\rangle_{2}\no\\
&&=\frac{1}{k}w({\cal O}_{e^{N-2-n}}{\cal O}_{e^{n-1-2(k-N)}}{\cal O}_{e})_{2}
-\frac{1}{k}\tilde{L}_{1+(k-N)}^{N,k,1}
w({\cal O}_{e^{N-2-n}}{\cal O}_{e^{n-1-2(k-N)}}{\cal O}_{e}
{\cal O}_{e^{1+(k-N)}})_{1}\no\\
&&=V_{2}^{N,k,2}(n;0)-\tilde{L}_{1+(k-N)}^{N,k,1}V_{1}^{N,k,2}(n;1)
-\tilde{L}_{1+(k-N)}^{N,k,1}V_{1}^{N,k,2}(n;1)\no\\
&&=V_{2}^{N,k,2}(n;0)-2\tilde{L}_{1+(k-N)}^{N,k,1}V_{1}^{N,k,2}(n;1).
\label{r2}
\end{eqnarray}
Then we turn into the $d=3$ case. In this case, we have to use the modified
K\"ahler equation (\ref{mok}), associativity equation and (\ref{wex2}).
\begin{eqnarray}
&&\frac{1}{k}\langle{\cal O}_{e^{N-2-n}}{\cal O}_{e^{n-1-3(k-N)}}{\cal O}_{e}
\rangle_{3}\no\\
&&=\frac{1}{k}\biggl(w({\cal O}_{e^{N-2-n}}{\cal O}_{e^{n-1-3(k-N)}}
{\cal O}_{e}
)_{3}-\tilde{L}^{N,k,1}_{1+(k-N)}
w({\cal O}_{e^{N-2-n}}{\cal O}_{e^{n-1-3(k-N)}}{\cal O}_{e}
{\cal O}_{e^{1+(k-N)}})_{2}\no\\
&&-\frac{1}{2}\tilde{L}^{N,k,2}_{1+2(k-N)}
w({\cal O}_{e^{N-2-n}}{\cal O}_{e^{n-1-3(k-N)}}{\cal O}_{e}
{\cal O}_{e^{1+2(k-N)}})_{1}\no\\
&&+\frac{1}{2}(\tilde{L}^{N,k,1}_{1+(k-N)})^{2}
w({\cal O}_{e^{N-2-n}}{\cal O}_{e^{n-1-3(k-N)}}{\cal O}_{e}
{\cal O}_{e^{1+(k-N)}}{\cal O}_{e^{1+(k-N)}})_{1}\biggr)\no\\
&&=\frac{1}{k}\biggl(w({\cal O}_{e^{N-2-n}}{\cal O}_{e^{n-1-3(k-N)}}
{\cal O}_{e})_{3}-2\tilde{L}^{N,k,1}_{1+(k-N)}
w({\cal O}_{e^{N-2-n}}{\cal O}_{e^{n-1-3(k-N)}}
{\cal O}_{e^{1+(k-N)}})_{2}\no\\
&&-\frac{1}{2}\tilde{L}^{N,k,2}_{1+2(k-N)}
w({\cal O}_{e^{N-2-n}}{\cal O}_{e^{n-1-3(k-N)}}
{\cal O}_{e^{1+2(k-N)}})_{1}\no\\
&&+\frac{3}{2}(\tilde{L}^{N,k,1}_{1+(k-N)})^{2}
w({\cal O}_{e^{N-2-n}}{\cal O}_{e^{n-1-3(k-N)}}
{\cal O}_{e^{1+(k-N)}}{\cal O}_{e^{1+(k-N)}})_{1}\biggr)\no\\
&&=V_{3}^{N,k,3}(n;0)-\tilde{L}^{N,k,1}_{1+(k-N)}V_{2}^{N,k,3}(n;1)
-\tilde{L}^{N,k,2}_{1+2(k-N)}V_{1}^{N,k,3}(n;2)
+(\tilde{L}^{N,k,1}_{1+(k-N)})^{2}V_{1}^{N,k,3}(n;1+1)\no\\
&&-2\tilde{L}^{N,k,1}_{1+(k-N)}\bigl(V_{2}^{N,k,3}(n;1)-
\tilde{L}^{N,k,1}_{1+(k-N)}V_{1}^{N,k,3}(n;1+1)\bigr)-
\frac{1}{2}\tilde{L}^{N,k,2}_{1+2(k-N)}V_{1}^{N,k,3}(n;2)\no\\
&&+\frac{3}{2}(\tilde{L}^{N,k,1}_{1+(k-N)})^{2}V_{1}^{N,k,3}(n;1+1)\no\\
&&=V_{3}^{N,k,3}(n;0)-3\tilde{L}^{N,k,1}_{1+(k-N)}V_{2}^{N,k,3}(n;1)
-\frac{3}{2}\tilde{L}^{N,k,2}_{1+2(k-N)}V_{1}^{N,k,3}(n;2)
+\frac{9}{2}(\tilde{L}^{N,k,1}_{1+(k-N)})^{2}V_{1}^{N,k,3}(n;1+1).\no\\
\end{eqnarray}
Indeed, these results agree with the results in \cite{vir}.
Next, we derive the generalized mirror transformation of $M_{k-1}^{k}$ 
model for the $d=4,5$ cases.
\begin{eqnarray}
&&\frac{1}{k}\langle{\cal O}_{e^{k-3-n}}{\cal O}_{e^{n-5}}
{\cal O}_{e}\rangle_{4}\no\\
&&=\frac{1}{k}\biggl(w({\cal O}_{e^{k-3-n}}{\cal O}_{e^{n-5}}
{\cal O}_{e})_{4}\nonumber\\
&&-\tilde{L}_{2}^{k-1,k,1}w({\cal O}_{e^{k-3-n}}{\cal O}_{e^{n-5}}
{\cal O}_{e}{\cal O}_{e^{2}})_{3}
-\frac{1}{2}\tilde{L}_{3}^{k-1,k,2}w({\cal O}_{e^{k-3-n}}{\cal O}_{e^{n-5}}
{\cal O}_{e}{\cal O}_{e^{3}})_{2}\nonumber\\
&&-\frac{1}{3}\tilde{L}_{4}^{k-1,k,3}w({\cal O}_{e^{k-3-n}}{\cal O}_{e^{n-5}}
{\cal O}_{e}{\cal O}_{e^{4}})_{1}\no\\
&&+\frac{1}{2}(\tilde{L}_{2}^{k-1,k,1})^{2}
w({\cal O}_{e^{k-3-n}}{\cal O}_{e^{n-5}}
{\cal O}_{e}{\cal O}_{e^{2}}{\cal O}_{e^{2}})_{2}
+\frac{1}{2}\tilde{L}_{2}^{k-1,k,1}\tilde{L}_{3}^{k-1,k,2}
w({\cal O}_{e^{k-3-n}}{\cal O}_{e^{n-5}}
{\cal O}_{e}{\cal O}_{e^{2}}{\cal O}_{e^{3}})_{1}\no\\
&&-\frac{1}{6}(\tilde{L}_{2}^{k-1,k,1})^{3}
w({\cal O}_{e^{k-3-n}}{\cal O}_{e^{n-6}}
{\cal O}_{e}{\cal O}_{e^{2}}{\cal O}_{e^{2}}{\cal O}_{e^{2}})_{1}\biggr)
\no\\
&&=\frac{1}{k}\biggl(w({\cal O}_{e^{k-3-n}}{\cal O}_{e^{n-5}}
{\cal O}_{e})_{4}\nonumber\\
&&-3\tilde{L}_{2}^{k-1,k,1}w({\cal O}_{e^{k-3-n}}{\cal O}_{e^{n-5}}
{\cal O}_{e^{2}})_{3}
-\tilde{L}_{3}^{k-1,k,2}w({\cal O}_{e^{k-3-n}}{\cal O}_{e^{n-5}}
{\cal O}_{e^{3}})_{2}\nonumber\\
&&-\frac{1}{3}\tilde{L}_{4}^{k-1,k,3}w({\cal O}_{e^{k-3-n}}{\cal O}_{e^{n-5}}
{\cal O}_{e^{4}})_{1}\no\\
&&+2(\tilde{L}_{2}^{k-1,k,1})^{2}
w({\cal O}_{e^{k-3-n}}{\cal O}_{e^{n-5}}
{\cal O}_{e^{2}}{\cal O}_{e^{2}})_{2}
+2\tilde{L}_{2}^{k-1,k,1}\tilde{L}_{3}^{k-1,k,2}
w({\cal O}_{e^{k-3-n}}{\cal O}_{e^{n-5}}
{\cal O}_{e^{2}}{\cal O}_{e^{3}})_{1}\no\\
&&-\frac{2}{3}(\tilde{L}_{2}^{k-1,k,1})^{3}
w({\cal O}_{e^{k-3-n}}{\cal O}_{e^{n-6}}
{\cal O}_{e^{2}}{\cal O}_{e^{2}}{\cal O}_{e^{2}})_{1}\biggr)\no\\
&&=V_{4}^{k-1,k,4}(n;0)-\tilde{L}^{k-1,k,1}_{2}V_{3}^{k-1,k,4}(n;1)
-\tilde{L}^{k-1,k,2}_{3}V_{2}^{k-1,k,4}(n;2)
-\tilde{L}^{k-1,k,3}_{4}V_{1}^{k-1,k,4}(n;3)\no\\
&&+(\tilde{L}^{k-1,k,1}_{2})^{2}\biggl(
V_{2}^{k-1,k,3}(n;1)+V_{2}^{k-1,k,3}(n-1;1)
-V_{2}^{k-1,k,2}(5;0)+A(n)\biggr)\no\\
&&+2\tilde{L}^{k-1,k,1}_{2}\tilde{L}^{k-1,k,2}_{3}
V_{1}^{k-1,k,4}(n;1+2)
-(\tilde{L}^{k-1,k,1}_{2})^{3}V_{1}^{k-1,k,4}(n;1+1+1)\no\\
&&-3\tilde{L}^{k-1,k,1}_{2}\biggl(V_{3}^{k-1,k,4}(n;1)-
\tilde{L}^{k-1,k,1}_{2}\bigl(V_{2}^{k-1,k,3}(n;1)+V_{2}^{k-1,k,3}(n-1;1)
-V_{2}^{k-1,k,2}(5;0)+A(n)\bigr)\no\\
&&-\tilde{L}^{k-1,k,2}_{3}V_{1}^{k-1,k,4}(n;1+2)
+(\tilde{L}^{k-1,k,1}_{2})^{2}V_{1}^{k-1,k,4}(n;1+1+1)\biggr)\no\\
&&-\tilde{L}^{k-1,k,2}_{3}\bigl(V_{2}^{k-1,k,4}(n;2)
-\tilde{L}^{k-1,k,1}_{2}V_{1}^{k-1,k,4}(n;1+2)\bigr)\no\\
&&-\frac{1}{3}\tilde{L}^{k-1,k,3}_{4}V_{1}^{k-1,k,4}(n;3)\no\\
&&+2(\tilde{L}^{k-1,k,1}_{2})^{2}
\bigl(2V_{2}^{k-1,k,3}(n;1)+2V_{2}^{k-1,k,3}(n-1;1)
-2V_{2}^{k-1,k,2}(5;0)+A(n)\no\\
&&-3\tilde{L}^{k-1,k,1}_{2}V_{1}^{k-1,k,4}(n;1+1+1)\bigr)
+2\tilde{L}^{k-1,k,1}_{2}\tilde{L}^{k-1,k,2}_{3}
V_{1}^{k-1,k,4}(n;1+2)\no\\
&&-\frac{2}{3}(\tilde{L}^{k-1,k,1}_{2})^{3}V_{1}^{k-1,k,4}(n;1+1+1)\no\\
&&=V_{4}^{k-1,k,4}(n;0)-4\tilde{L}^{k-1,k,1}_{2}V_{3}^{k-1,k,4}(n;1)
-2\tilde{L}^{k-1,k,2}_{3}V_{2}^{k-1,k,4}(n;2)
-\frac{4}{3}\tilde{L}^{k-1,k,3}_{4}V_{1}^{k-1,k,4}(n;3)\no\\
&&+8(\tilde{L}^{k-1,k,1}_{2})^{2}\biggl(
V_{2}^{k-1,k,3}(n;1)+V_{2}^{k-1,k,3}(n-1;1)
-V_{2}^{k-1,k,2}(5;0)+\frac{3}{4}A(n)\biggr)\no\\
&&+8\tilde{L}^{k-1,k,1}_{2}\tilde{L}^{k-1,k,2}_{3}
V_{1}^{k-1,k,4}(n;1+2)
-\frac{32}{3}(\tilde{L}^{k-1,k,1}_{2})^{3}V_{1}^{k-1,k,4}(n;1+1+1)
\label{deri4}
\end{eqnarray}
Note that we have used (\ref{mok}) to derive the second equality in 
(\ref{deri4}).
In this derivation, the mysterious change of coefficients $\frac{1}{2}$
of $A(n)$ in (\ref{virc}) into $\frac{3}{4}$ in (\ref{virg}) naturally 
arises as the result of the identity:
\begin{equation}
1+3+4\cdot\frac{1}{2}=\frac{4^{2}}{2!}\cdot\frac{3}{4}.
\label{mys1}
\end{equation} 
Derivation of the $d=5$ result can be done in almost the same way 
as in the $d=4$ case. As the first step, we use (\ref{mok}) once again,
\begin{eqnarray}
&&\frac{1}{k}\langle{\cal O}_{e^{k-3-n}}{\cal O}_{e^{n-6}}
{\cal O}_{e}\rangle_{5}\no\\
&&=\frac{1}{k}\biggl(w({\cal O}_{e^{k-3-n}}{\cal O}_{e^{n-6}}
{\cal O}_{e})_{5}\nonumber\\
&&-\tilde{L}_{2}^{k-1,k,1}w({\cal O}_{e^{k-3-n}}{\cal O}_{e^{n-6}}
{\cal O}_{e}{\cal O}_{e^{2}})_{4}
-\frac{1}{2}\tilde{L}_{3}^{k-1,k,2}w({\cal O}_{e^{k-3-n}}{\cal O}_{e^{n-6}}
{\cal O}_{e}{\cal O}_{e^{3}})_{3}\nonumber\\
&&-\frac{1}{3}\tilde{L}_{4}^{k-1,k,3}w({\cal O}_{e^{k-3-n}}{\cal O}_{e^{n-6}}
{\cal O}_{e}{\cal O}_{e^{4}})_{2}
-\frac{1}{4}\tilde{L}_{5}^{k-1,k,4}w({\cal O}_{e^{k-3-n}}{\cal O}_{e^{n-6}}
{\cal O}_{e}{\cal O}_{e^{5}})_{1}\nonumber\\
&&+\frac{1}{2}(\tilde{L}_{2}^{k-1,k,1})^{2}
w({\cal O}_{e^{k-3-n}}{\cal O}_{e^{n-6}}
{\cal O}_{e}{\cal O}_{e^{2}}{\cal O}_{e^{2}})_{3}
+\frac{1}{2}\tilde{L}_{2}^{k-1,k,1}\tilde{L}_{3}^{k-1,k,2}
w({\cal O}_{e^{k-3-n}}{\cal O}_{e^{n-6}}
{\cal O}_{e}{\cal O}_{e^{2}}{\cal O}_{e^{3}})_{2}\no\\
&&+\frac{1}{3}\tilde{L}_{2}^{k-1,k,1}\tilde{L}_{4}^{k-1,k,3}
w({\cal O}_{e^{k-3-n}}{\cal O}_{e^{n-6}}
{\cal O}_{e}{\cal O}_{e^{2}}{\cal O}_{e^{4}})_{1}
+\frac{1}{8}(\tilde{L}_{3}^{k-1,k,2})^{2}
w({\cal O}_{e^{k-3-n}}{\cal O}_{e^{n-6}}
{\cal O}_{e}{\cal O}_{e^{3}}{\cal O}_{e^{3}})_{1}\nonumber\\
&&-\frac{1}{6}(\tilde{L}_{2}^{k-1,k,1})^{3}
w({\cal O}_{e^{k-3-n}}{\cal O}_{e^{n-6}}
{\cal O}_{e}{\cal O}_{e^{2}}{\cal O}_{e^{2}}{\cal O}_{e^{2}})_{2}
-\frac{1}{4}(\tilde{L}_{2}^{k-1,k,1})^{2}\tilde{L}_{3}^{k-1,k,2}
w({\cal O}_{e^{k-3-n}}{\cal O}_{e^{n-6}}
{\cal O}_{e}{\cal O}_{e^{2}}{\cal O}_{e^{2}}{\cal O}_{e^{3}})_{1}\nonumber\\
&&+\frac{1}{24}(\tilde{L}_{2}^{k-1,k,1})^{4}
w({\cal O}_{e^{k-3-n}}{\cal O}_{e^{n-6}}
{\cal O}_{e}{\cal O}_{e^{2}}{\cal O}_{e^{2}}{\cal O}_{e^{2}}{\cal O}_{e^{2}}
)_{1}\biggr)\nonumber\\
&&=\frac{1}{k}\biggl(w({\cal O}_{e^{k-3-n}}{\cal O}_{e^{n-6}}
{\cal O}_{e})_{5}\nonumber\\
&&-4\tilde{L}_{2}^{k-1,k,1}w({\cal O}_{e^{k-3-n}}{\cal O}_{e^{n-6}}
{\cal O}_{e^{2}})_{4}
-\frac{3}{2}\tilde{L}_{3}^{k-1,k,2}w({\cal O}_{e^{k-3-n}}{\cal O}_{e^{n-6}}
{\cal O}_{e^{3}})_{3}\nonumber\\
&&-\frac{2}{3}\tilde{L}_{4}^{k-1,k,3}w({\cal O}_{e^{k-3-n}}{\cal O}_{e^{n-6}}
{\cal O}_{e^{4}})_{2}
-\frac{1}{4}\tilde{L}_{5}^{k-1,k,4}w({\cal O}_{e^{k-3-n}}{\cal O}_{e^{n-6}}
{\cal O}_{e^{5}})_{1}\nonumber\\
&&+\frac{5}{2}(\tilde{L}_{2}^{k-1,k,1})^{2}
w({\cal O}_{e^{k-3-n}}{\cal O}_{e^{n-6}}
{\cal O}_{e^{2}}{\cal O}_{e^{2}})_{3}
+\frac{5}{2}\tilde{L}_{2}^{k-1,k,1}\tilde{L}_{3}^{k-1,k,2}
w({\cal O}_{e^{k-3-n}}{\cal O}_{e^{n-6}}
{\cal O}_{e^{2}}{\cal O}_{e^{3}})_{2}\nonumber\\
&&+\frac{5}{3}\tilde{L}_{2}^{k-1,k,1}\tilde{L}_{4}^{k-1,k,3}
w({\cal O}_{e^{k-3-n}}{\cal O}_{e^{n-6}}
{\cal O}_{e^{2}}{\cal O}_{e^{4}})_{1}
+\frac{5}{8}(\tilde{L}_{3}^{k-1,k,2})^{2}
w({\cal O}_{e^{k-3-n}}{\cal O}_{e^{n-6}}
{\cal O}_{e^{3}}{\cal O}_{e^{3}})_{3}\nonumber\\
&&-\frac{5}{6}(\tilde{L}_{2}^{k-1,k,1})^{3}
w({\cal O}_{e^{k-3-n}}{\cal O}_{e^{n-6}}
{\cal O}_{e^{2}}{\cal O}_{e^{2}}{\cal O}_{e^{2}})_{2}
-\frac{5}{4}(\tilde{L}_{2}^{k-1,k,1})^{2}\tilde{L}_{3}^{k-1,k,2}
w({\cal O}_{e^{k-3-n}}{\cal O}_{e^{n-6}}
{\cal O}_{e^{2}}{\cal O}_{e^{2}}{\cal O}_{e^{3}})_{1}\nonumber\\
&&+\frac{5}{24}(\tilde{L}_{2}^{k-1,k,1})^{4}
w({\cal O}_{e^{k-3-n}}{\cal O}_{e^{n-6}}
{\cal O}_{e^{2}}{\cal O}_{e^{2}}{\cal O}_{e^{2}}{\cal O}_{e^{2}}
)_{1}\biggr).\nonumber\\
\label{deri5}
\end{eqnarray}
Then we decompose each $w({\cal O}_{e^{n_{1}}}{\cal O}_{e^{n_{2}}}
\cdots {\cal O}_{e^{n_{l}}})$ in terms of 
$V_{d-m}^{k-1,k,d}(n;\sigma_{m})$ by using the associativity equation 
and the modified K\"ahler equation (\ref{mok}). We briefly explain the 
remaining computation. At first, the coefficients of the terms in 
(\ref{deri5}) add up to the ones in (\ref{gene}) as follows,
\begin{eqnarray}
&&1+4=5,\;\;1+\frac{3}{2}=\frac{5}{2},\;\;1+\frac{2}{3}=\frac{5}{3},\;\;
1+\frac{1}{4}=\frac{5}{4},\no\\
&&1+4+3\cdot\frac{5}{2}=\frac{25}{2},\;\;2+4+\frac{3}{2}+\frac{5}{2}\cdot2=
\frac{25}{2},\;\;2+4+\frac{2}{3}+\frac{5}{3}=\frac{25}{3},\;\;
1+\frac{3}{2}+\frac{5}{8}=\frac{25}{8},\no\\
&&1+4+3\cdot\frac{5}{2}+2\cdot\frac{5}{2}+2^{2}\cdot
\frac{5}{6}=\frac{125}{6},\;\;
3+2\cdot4+\frac{3}{2}+3\cdot\frac{5}{2}+2\cdot\frac{5}{2}+\frac{5}{2}+
\frac{5}{2}+\frac{5}{4}=\frac{125}{4},\no\\
&&1+4+3\cdot\frac{5}{2}+2\cdot\frac{5}{2}+2^{2}\cdot\frac{5}{6}+\frac{5}{2}
+2\cdot\frac{5}{6}+\frac{5}{6}+\frac{5}{24}=\frac{625}{24}.
\end{eqnarray}
The change of coefficients $\frac{2}{3},
\frac{1}{3},\frac{1}{2}$ of $B(n)$, $C(n)$, $D(n)$ in 
(\ref{virc}) into $\frac{4}{5},
\frac{3}{5},\frac{4}{5}$ in (\ref{virg}) 
is naturally derived in the same way as in the 
$d=4$ case:
\begin{eqnarray}
&&1+4+3\cdot\frac{5}{2}\cdot\frac{2}{3}=\frac{25}{2}\cdot\frac{4}{5},
\;\;1+4+3\cdot\frac{5}{2}\cdot\frac{1}{3}=\frac{25}{2}\cdot\frac{3}{5},\no\\
&&1+4+3\cdot\frac{5}{2}+(2\cdot\frac{5}{2}+2^{2}\cdot
\frac{5}{6})\cdot\frac{1}{2}=\frac{125}{6}\cdot\frac{4}{5}.
\end{eqnarray}
Lastly, we verify the non-trivial change of the terms including 
$hi_{1}(n),hi_{2}(n),hi_{3}(n),hi_{4}(n)$ in (\ref{virg}). 
We discuss the case of $1+1+1$ sector in detail. In this case, 
we have four terms in (\ref{deri5}) that contribute to this sector:
\begin{eqnarray}
&&\frac{1}{k}\cdot w({\cal O}_{e^{k-3-n}}{\cal O}_{e^{n-6}}
{\cal O}_{e})_{5},\;\;
-\frac{1}{k}\cdot4\tilde{L}_{2}^{k-1,k,1}w({\cal O}_{e^{k-3-n}}{\cal O}_{e^{n-6}}
{\cal O}_{e^{2}})_{4}\no\\
&&\frac{1}{k}\cdot\frac{5}{2}(\tilde{L}_{2}^{k-1,k,1})^{2}
w({\cal O}_{e^{k-3-n}}{\cal O}_{e^{n-6}}
{\cal O}_{e^{2}}{\cal O}_{e^{2}})_{3},\;\;
-\frac{1}{k}\cdot\frac{5}{6}(\tilde{L}_{2}^{k-1,k,1})^{3}
w({\cal O}_{e^{k-3-n}}{\cal O}_{e^{n-6}}
{\cal O}_{e^{2}}{\cal O}_{e^{2}}{\cal O}_{e^{2}})_{2}
\end{eqnarray} 
If we consider the contributions coming from the first two terms, we can see 
that we don't have to use the modified K\"ahler equation (\ref{mok}). 
Therefore, we can compute the contributions that includes $hi_{j}(n)$ only 
by using associativity equation and (\ref{wex}). The contributions from 
these terms turn out to be the same and given by,
\begin{eqnarray}
&&V_{1}^{k-1,k,3}(n;1+1)\cdot V_{1}^{k-1,k,1}(n-4;0)
+V_{1}^{k-1,k,2}(n;1)\cdot V_{1}^{k-1,k,2}(n-3;1)\no\\
&&+V_{1}^{k-1,k,1}(n;0)\cdot V_{1}^{k-1,k,3}(n-2;1+1)\no\\
&&-V_{1}^{k-1,k,5}(n;1+1+2)\cdot V_{1}^{k-1,k,1}(3;0)
-V_{1}^{k-1,k,5}(n;1+3)\cdot V_{1}^{k-1,k,2}(4;1)\no\\
&&-V_{1}^{k-1,k,5}(n;4)\cdot V_{1}^{k-1,k,3}(5;1+1)\no\\
&&=3hi_{1}(n)+3hi_{2}(n)+hi_{3}(n).
\label{dhi}
\end{eqnarray}
Contrary to the above cases, we have to treat carefully the term  
$w({\cal O}_{e^{k-3-n}}{\cal O}_{e^{n-6}}
{\cal O}_{e^{2}}{\cal O}_{e^{2}})_{3}$ because we have to use (\ref{mok}) once 
in the computation. In this case, we take care of the following terms that 
appear in decomposing $w({\cal O}_{e^{k-3-n}}{\cal O}_{e^{n-6}}
{\cal O}_{e^{2}}{\cal O}_{e^{2}})_{3}$ by using associativity 
equation:
\begin{eqnarray}
&&\frac{1}{k^{2}}\biggl(w({\cal O}_{e^{k-3-n}}{\cal O}_{e^{2}}
{\cal O}_{e}{\cal O}_{e^{n-4}})_{2}\cdot 
w({\cal O}_{e^{k-n+1}}{\cal O}_{e^{n-6}}{\cal O}_{e})_{1}
+w({\cal O}_{e^{k-3-n}}{\cal O}_{e^{2}}
{\cal O}_{e}{\cal O}_{e^{n-3}})_{1}\cdot 
w({\cal O}_{e^{k-n}}{\cal O}_{e^{n-6}}{\cal O}_{e})_{2}\no\\
&&+w({\cal O}_{e^{k-3-n}}
{\cal O}_{e}{\cal O}_{e^{n-3}})_{2}\cdot 
w({\cal O}_{e^{k-n}}{\cal O}_{e^{n-6}}{\cal O}_{e^{2}}{\cal O}_{e})_{1}
+w({\cal O}_{e^{k-3-n}}
{\cal O}_{e}{\cal O}_{e^{n-2}})_{1}\cdot 
w({\cal O}_{e^{k-n-1}}{\cal O}_{e^{n-6}}{\cal O}_{e^{2}}{\cal O}_{e})_{2}\no\\
&&-w({\cal O}_{e^{k-3-n}}
{\cal O}_{e^{n-6}}{\cal O}_{e^{2}}{\cal O}_{e^{3}} )_{2}\cdot 
w({\cal O}_{e^{k-6}}{\cal O}_{e}{\cal O}_{e})_{1}
-w({\cal O}_{e^{k-3-n}}
{\cal O}_{e^{n-6}}{\cal O}_{e^{2}}{\cal O}_{e^{4}} )_{1}\cdot 
w({\cal O}_{e^{k-7}}{\cal O}_{e}{\cal O}_{e})_{2}\no\\
&&-w({\cal O}_{e^{k-3-n}}
{\cal O}_{e^{n-6}}{\cal O}_{e^{4}} )_{2}\cdot 
w({\cal O}_{e^{k-7}}{\cal O}_{e^{2}}{\cal O}_{e}{\cal O}_{e})_{1}
-w({\cal O}_{e^{k-3-n}}
{\cal O}_{e^{n-6}}{\cal O}_{e^{5}} )_{1}\cdot 
w({\cal O}_{e^{k-8}}{\cal O}_{e^{2}}{\cal O}_{e}{\cal O}_{e})_{1}\biggr).
\label{A}
\end{eqnarray}
If we take the leading terms of (\ref{mok}) in reducing the degree $2$ 
four point 
functions in (\ref{A}) into three point functions and if we 
pick up the terms proportional to $\tilde{L}^{k-1,k,1}_{2}$ 
in these degree $2$ three point function, we obtain the 
following contribution:
\begin{eqnarray}
2(3hi_{1}(n)+3hi_{2}(n)+hi_{3}(n)).
\label{ddhi}
\end{eqnarray} 
On the other hand, we also have the terms proportional to 
$\tilde{L}^{k-1,k,1}_{2}$ by picking up the sub-leading term of (\ref{mok}) 
that appear in reducing  the degree $2$ four point function in (\ref{A}).
For example, we take \\
$\tilde{L}^{k-1,k,1}_{2}w({\cal O}_{e^{k-1-n}}
{\cal O}_{e^{2}}{\cal O}_{e^{2}}{\cal O}_{e^{n-4}} )_{1}$ in the r.h.s. of 
 the equation:
\begin{equation}
w({\cal O}_{e^{k-1-n}}
{\cal O}_{e^{2}}{\cal O}_{e}{\cal O}_{e^{n-4}} )_{2}=
2w({\cal O}_{e^{k-1-n}}
{\cal O}_{e^{2}}{\cal O}_{e^{n-4}} )_{2}
-\tilde{L}^{k-1,k,1}_{2}
w({\cal O}_{e^{k-1-n}}
{\cal O}_{e^{2}}{\cal O}_{e^{2}}{\cal O}_{e^{n-4}} )_{1}.
\label{emok}
\end{equation}
After some computation, the following contribution appears,
\begin{eqnarray}
&&V_{1}^{k-1,k,3}(n;1+1)\cdot V_{1}^{k-1,k,1}(n-4;0)
+V_{1}^{k-1,k,1}(n;0)\cdot V_{1}^{k-1,k,3}(n-2;1+1)\no\\
&&-V_{1}^{k-1,k,5}(n;1+1+2)\cdot V_{1}^{k-1,k,1}(3;0)
-V_{1}^{k-1,k,5}(n;4)\cdot V_{1}^{k-1,k,3}(5;1+1)\no\\
&&=2hi_{1}(n)+2hi_{2}(n)+hi_{3}(n)-hi_{4}(n).
\label{ghi}
\end{eqnarray}
The contribution from $w({\cal O}_{e^{k-3-n}}{\cal O}_{e^{n-6}}
{\cal O}_{e^{2}}{\cal O}_{e^{2}}{\cal O}_{e^{2}})_{2}$ is obtained 
only by taking leading terms of (\ref{mok}) and of three 
point functions. Therefore, the computation is the same as the one of virtual 
Gromov-Witten invariants, and we have from (\ref{virc}),
\begin{eqnarray}
&&V_{1}^{k-1,k,3}(n;1+1)\cdot V_{1}^{k-1,k,1}(n-4;0)
+2V_{1}^{k-1,k,2}(n;1)\cdot V_{1}^{k-1,k,2}(n-3;1)\no\\
&&+V_{1}^{k-1,k,1}(n;0)\cdot V_{1}^{k-1,k,3}(n-2;1+1)\no\\
&&-V_{1}^{k-1,k,5}(n;1+1+2)\cdot V_{1}^{k-1,k,1}(3;0)
-2V_{1}^{k-1,k,5}(n;1+3)\cdot V_{1}^{k-1,k,2}(4;1)\no\\
&&-V_{1}^{k-1,k,5}(n;4)\cdot V_{1}^{k-1,k,3}(5;1+1)\no\\
&&=4hi_{1}(n)+4hi_{2}(n)+hi_{3}(n)+hi_{4}(n).
\label{hhi}
\end{eqnarray}  
Adding up these five contributions, we can finally reproduce the part of  
the formula of
$G_{2}^{k-1,k,5}(n;1+1+1)$ including $hi_{j}(n)$:     
\begin{eqnarray}
&&(3hi_{1}(n)+3hi_{2}(n)+hi_{3}(n))
+4(3hi_{1}(n)+3hi_{2}(n)+hi_{3}(n))\no\\
&&+\frac{5}{2}\cdot2(3hi_{1}(n)+3hi_{2}(n)+hi_{3}(n))
+\frac{5}{2}
(2hi_{1}(n)+2hi_{2}(n)+hi_{3}(n)-hi_{4}(n))\no\\
&&+\frac{5}{6}(4hi_{1}(n)+4hi_{2}(n)+hi_{3}(n)+hi_{4}(n))\no\\
&&=\frac{125}{6}
(\frac{46}{25}hi_{1}(n)+\frac{46}{25}hi_{2}(n)+\frac{16}{25}hi_{3}(n)-
\frac{2}{25}hi_{4}(n)).
\label{fin}
\end{eqnarray}
Computation of the remaining $1+1$, $1+2$ sectors goes in the same way. 
In these cases, we don't have to consider the sub-leading terms of (\ref{mok})
because we don't have to use (\ref{mok}) except for 
$w({\cal O}_{e^{k-3-n}}{\cal O}_{e^{n-6}}
{\cal O}_{e^{2}}{\cal O}_{e^{2}})_{3}$ and 
$w({\cal O}_{e^{k-3-n}}{\cal O}_{e^{n-6}}
{\cal O}_{e^{2}}{\cal O}_{e^{3}})_{2}$. Therefore, we write down the final 
process of the computation in the following:
\begin{eqnarray}
&&1+2\;\;\mbox{sector}\no\\
&&2(2hi_{1}(n)+\frac{3}{2}hi_{2}(n)+hi_{3}(n)-\frac{1}{2}hi_{4}(n))
+4(2hi_{1}(n)+hi_{2}(n)+hi_{3}(n)-hi_{4}(n))\no\\
&&+\frac{3}{2}(2hi_{1}(n)+2hi_{2}(n)+hi_{3}(n))
+\frac{5}{2}(2hi_{1}(n)+hi_{2}(n)+hi_{3}(n)-hi_{4}(n))
\no\\
&&=\frac{25}{2}
(\frac{8}{5}hi_{1}(n)+hi_{2}(n)+\frac{4}{5}hi_{3}(n)-\frac{3}{5}hi_{4}(n)),
\no\\
&&1+1\;\;\mbox{sector}\no\\
&&(2hi_{1}(n)+2hi_{2}(n)+hi_{3}(n))+4(2hi_{1}(n)+2hi_{2}(n)+hi_{3}(n))\no\\
&&+\frac{5}{2}\cdot(2hi_{1}(n)+hi_{2}(n)+hi_{3}(n)-hi_{4}(n))\no\\
&&=\frac{25}{2}
(\frac{6}{5}hi_{1}(n)+hi_{2}(n)+\frac{3}{5}hi_{3}(n)-\frac{1}{5}hi_{4}(n)).
\end{eqnarray}
In this way, we have reproduced the results in \cite{vir}.
We also did some numerical test for curves of higher degree.
In the $d=6$ case, Conjecture 3 predicts 
\begin{eqnarray}
&&L_{8}^{13,14,6}=
389591981138964503377056394266126437146547452695642109769122604114\setminus
\no\\
&&0067266620858637887736545388962432/9375,
\end{eqnarray}
which agrees with the numerical computation using the fixed-point theorem
\cite{kont}.
\section{Iritani's Theory}
In this section, we briefly explain how 
our conjecture (\ref{main}) naturally follows from applying Iritani's result 
\cite{iri} to our specific case. 
Crucial point of his framework is introduction of deformation variable 
of $x^{j}\;(j=1,\cdots, N-2),\;\;(x^{1}=x)$ that corresponds to the 
insertion of ${\cal O}_{e^{j}}$. 
First, we prepare a $(N-1)\times (N-1)$ 
matrix $\tilde{C}_{1}(e^{x^{1}})$ whose matrix element is given by the 
virtual structure constants:
\begin{equation}
(\tilde{C}_{1}(e^{x^{1}}))_{m}{}^{n}=\left\{ \begin{array}{ll}
\tilde{L}_{m}^{N,k,d}e^{dx^{1}}&\;\;(n=m+1+(k-N)d)\\
0&\;\;\mbox{otherwise}.\end{array}\right.
\label{c1}
\end{equation}
 Then we can express the truncated Gauss-Manin system (\ref{tgm}) into 
a compact form:
\begin{equation}
\d_{x^{1}}\vec{\psi}
=\tilde{C}_{1}(e^{x^{1}})\vec{\psi}.
\label{ctgm0}
\end{equation}
 Next, we construct matrices $\tilde{C}_{j}(e^{x^{1}})\;
(j=1,\cdots,N-2)$ that correspond to the deformation of $\vec{\psi}$
by $x^{j}$:
\begin{equation}
\d_{x^{j}}\vec{\psi}
=\tilde{C}_{j}(e^{x^{1}})\vec{\psi},
\label{ctgm}
\end{equation}
whose $mn$ element vanishes unless $n=m+j+(k-N)d$. Non-zero matrix elements 
of $C_{j}(e^{x^{1}})\;(j=1,\cdots,N-2)$ are uniquely determined 
by the conditions:
\begin{eqnarray}
[\tilde{C}_{i}(e^{x^{1}}),\tilde{C}_{j}(e^{x^{1}})]=0,\;\;
(\tilde{C}_{j}(e^{x^{1}}))_{m}{}^{m+j}=1.
\label{ini}
\end{eqnarray} 
With these set up, we can obtain Jacobi matrix between B-model deformation 
parameters $x^{j}$ and flat coordinates $t^{j}$ from the flat coordinate 
condition ($\frac{\d\psi_{0}}{\d t^{j}}=\psi_{j}$):
\begin{eqnarray}
&&\frac{\d\psi_{0}}{\d x^{i}}=\frac{\d t^{j}}{\d x^{i}}
\frac{\d\psi_{0}}{\d t^{j}}=\frac{\d t^{j}}{\d x^{i}}
\psi_{j}=(\tilde{C}_{i}(e^{x^{1}}))_{0}{}^{j}\psi_{j}\no\\
&&\Longrightarrow \frac{\d t^{j}}{\d x^{i}}=(\tilde{C}_{i}(e^{x^{1}}))_{0}
{}^{j}.
\label{iflat}
\end{eqnarray}
Using (\ref{iflat}), we can compute matrix $\bar{C}_{i}(e^{x^{1}})$ 
that corresponds to deformation of $t^{i}$ 
as follows:
\begin{equation}
\d_{t^{i}}\vec{\psi}=\frac{\d x^{j}}{\d t^{i}}\d_{x^{j}}\vec{\psi}
=\frac{\d x^{j}}{\d t^{i}}\tilde{C}_{j}(e^{x^{1}})\vec{\psi}
=\bar{C}_{i}(e^{x^{1}})\vec{\psi}
\label{fmat}
\end{equation}
Of course, we can construct rank 3 symmetric tensor 
\begin{equation}
\bar{C}_{ijm}(e^{x^{1}})=(\bar{C}_{i}(e^{x^{1}}))_{j}{}^{l}\eta_{lm},\;\;
(\eta_{lm}=k\cdot\delta_{l+m,N-2}).
\label{sym3}
\end{equation} 
Surprisingly, the non-zero elements of $\bar{C}_{1jm}(e^{x^{1}})$
is given by 
\begin{equation}
\bar{C}_{1,N-2-n,n-1-(k-N)d}(e^{x^{1}})
=w({\cal O}_{e^{N-2-n}}{\cal O}_{e^{n-1-(k-N)d}}
{\cal O}_{e})_{d},
\end{equation}
that was obtained by the first step of our construction.
On the other hand, if we integrate out the Jacobi matrix (\ref{iflat}), 
we can observe that expansion point of $x$-coordinates $(e^{x^{1}},0,\cdots,0)$
is given in terms of  $t$-coordinates as follows:
\begin{eqnarray}
&&(e^{x^{1}},0,\cdots,0)=(e^{t^{1}},t^{2}(t^{1}),t^{3}(t^{1}),\cdots,
t^{N-2}(t^{1})),\no\\
&&t^{i}(t^{1})=0,\;\;(i\notin\{1+(k-N)j\;\;|\;\;j\in {\bf N}\;\;\}),\;\;
t^{1+(k-N)j}(t^{1})=\frac{\tilde{L}_{1+(k-N)j}^{N,k,j}}{j}\exp(jt^{1})
,\;\;(j\in {\bf N}\;\;).
\label{itrans} 
\end{eqnarray}
Therefore, we have to change expansion point of $\bar{C}_{ijm}$ , i.e., we 
assert the following equation:
\begin{eqnarray}
&&C_{ijm}(e^{t^{1}},0,\cdots,0)=
\bar{C}_{ijm}(e^{t^{1}},-t^{2}(t^{1}),\cdots,-t^{N-2}(t^{1})).
\label{imt}
\end{eqnarray}
where $C_{ijm}(e^{t^{1}},0,\cdots,0)$ 
is the true three point function we want.
At this stage, we can easily see that the second step of our construction 
is nothing but the 
process of perturbing $\bar{C}_{ijm}(e^{x^{1}})$ by $x_{2},x_{3},\cdots,
x_{N-2}$ by using the associativity equation:
\begin{eqnarray}
&&\bar{C}_{ijm}(e^{x^{1}},x_{2},x_{3},\cdots,
x_{N-2})\eta^{ml}\bar{C}_{lst}(e^{x^{1}},x_{2},x_{3},\cdots,
x_{N-2})\no\\
&&=\bar{C}_{ism}(e^{x^{1}},x_{2},x_{3},\cdots,
x_{N-2})\eta^{ml}\bar{C}_{ljt}(e^{x^{1}},x_{2},x_{3},\cdots,
x_{N-2}),\no\\
&&\eta^{ml}=\frac{1}{k}\cdot\delta_{l+m,N-2},
\label{iaso}
\end{eqnarray}
and the K\"ahler equation: 
\begin{eqnarray}
&&\frac{d}{d t^{1}}\bar{C}_{ijm}(t^{1},-t^{2}(t_{1}),\cdots,-t^{N-2}(t_{1}))
\no\\
&&=\frac{\d}{\d x^{1}}
\bar{C}_{ijm}(e^{x^{1}},x^{2},\cdots,x^{N-2})
-\sum_{d=1}^{\infty}\tilde{L}_{1+(k-N)d}^{N,k,d}\exp(dt^{1})
\frac{\d}{\d x^{1+(k-N)d}}\bar{C}_{ijm}(e^{x^{1}},x^{2},\cdots,x^{N-2}).\no\\
\label{dmok}
\end{eqnarray}
Note that (\ref{dmok}) is nothing but the modified K\"ahler 
equation (\ref{mok}).
Lastly, we expand the r.h.s. of (\ref{imt}) in terms of $x^{2},\cdots,x^{N-2}$
and obtain (\ref{main}) as the 
special case of $i=1,\;j=N-2-n,\;m=n-1-(k-N)d$. 


\newpage
\begin{thebibliography}{99}
\bibitem{bert} A.Bertram.
\newblock{\em Another way to enumerate rational curves with torus
 actions}
\newblock{Invent.Math. {\bf 142} (2000), no.3, 487-512.}
\bibitem{givc} Tom Coates, Alexander B. Givental.
\newblock{\em Quantum Riemann-Roch, Lefschetz and Serre}
\newblock{math.AG/0110142}
\bibitem{cj}A. Collino, M.Jinzenji.
\newblock{\em On the Structure of Small Quantum Cohomology Rings for 
Projective Hypersurfaces}
\newblock{Commun.Math.Phys.206:157-183,1999}
\bibitem{gath} A. Gathmann.
{\em  Absolute and Relative Gromov-Witten Invariants of Very Ample 
Hypersurfaces}
\newblock math.AG/9908054.
\bibitem{giv}Alexander B. Givental.
\newblock{\em Equivariant Gromov - Witten Invariants}
\newblock{Internat. Math. Res.Notices 13 (1996),613--663.}
\bibitem{iri} H.Iritani. 
\newblock{\em Quantum D-module and Generalized Mirror Transformation }
\newblock{in preparation.}
\bibitem{fano} M.Jinzenji. 
\newblock{\em Completion of the Conjecture: 
Quantum Cohomology of Fano Hypersurfaces}
\newblock{Mod.Phys.Lett. A15 (2000) 101-120.}
\bibitem{vir}M.Jinzenji.
\newblock{\em Virtual Gromov-Witten Invariants and the Quantum Cohomology 
Rings of General Type Projective Hypersurfaces}
\newblock{Mod.Phys.Lett. A15 (2000) 629-650.}
\bibitem{gene}M.Jinzenji. 
\newblock{\em On the Quantum Cohomology Rings of General Type Projective 
Hypersurfaces and Generalized Mirror Transformation}
\newblock{Int.J.Mod.Phys.A15:1557-1596,2000.} 
\bibitem{gm}M.Jinzenji.
\newblock{\em Gauss-Manin System and the Virtual Structure Constants}
\newblock{Int.J.Math. 13 (2002) 445-478.}
\bibitem{kont} M.Kontsevich.
\newblock{\em Enumeration of Rational Curves via Torus Actions}
\newblock{The moduli space of curves, R.Dijkgraaf, C.Faber, G.van der
 Geer (Eds.), Progress in Math., v.129, Birkh\"auser, 1995, 335-368.}
\bibitem{km} M. Kontsevich, Yu. Manin.
\newblock{\em Gromov-Witten classes, quantum cohomology, 
and enumerative geometry}
\newblock{Commun.Math.Phys. 164 (1994) 525-562.}
\bibitem{blly}B.Lian, K.Liu and S.T.Yau.
\newblock{\em Mirror Principle III}
\newblock{Asian J. Math. {\bf 3} (1999), no.4, 771-800.}
\end{thebibliography}
\end{document}